\documentclass{article}
\usepackage{euscript,amsfonts,amssymb,amsmath,amscd,amsthm,url,graphicx,epsfig,ccaption}
\usepackage[all]{xy}
\author{Alexey G. Gorinov}
\title{Rational cohomology of the moduli spaces of pointed genus 1 curves}
\date{}
\newcommand{\eu}{\EuScript}
\newcommand{\Z}{{\mathbb Z}}
\newcommand{\C}{{\mathbb C}}
\newcommand{\R}{{\mathbb R}}
\newcommand{\Q}{{\mathbb Q}}

\newcommand{\aaa}{\mathbf{a}}
\newcommand{\bb}{\mathbf{b}}

\newcommand{\ee}{\mathbf{e}}

\newcommand{\Id}{\mathrm{Id}}

\newtheorem{theorem}{Theorem}
\newtheorem{lemma}{Lemma}
\newtheorem{Prop}{Proposition}

\newtheorem*{observation}{Observation}
\newtheorem{corollary}{Corollary}
\DeclareMathOperator\res{res}
\DeclareMathOperator\tr{tr}
\DeclareMathOperator\coker{coker}
\DeclareMathOperator\Gr{Gr}
\DeclareMathOperator\SL{SL}
\DeclareMathOperator\Imm{Im}
\DeclareMathOperator\PSL{PSL}
\DeclareMathOperator\Map{Map}
\DeclareMathOperator\Hom{Hom}
\DeclareMathOperator\Sym{Sym}
\DeclareMathOperator\Aut{Aut}
\DeclareMathOperator\End{End}

\begin{document}
\maketitle
\begin{abstract}
We give a combinatorial description of the rational cohomology of the moduli spaces of pointed genus 1 curves with $n$ marked points and level $N$ structures. More precisely, we explicitly describe the $E_2$ term of the Leray spectral sequence of the forgetful mapping $\EuScript{M}_{1,n}(N)\to\EuScript{M}_{1,1}(N)$ and show that the result is isomorphic to the rational cohomology of $\EuScript{M}_{1,n}(N)$ as a rational mixed Hodge structure equipped with an action of the symmetric group $\mathfrak{S}_n$. The classical moduli space $\EuScript{M}_{1,n}$ is the particular case $N=1$.
\end{abstract}
\section{Introduction}
We denote by $\eu{M}_{g,n}(N)$ the coarse moduli space of smooth genus $g$ complex curves with $n$ marked (ordered) points and level $N$ structures.
The purpose of this note is to reduce the computation of the rational cohomology of $\eu{M}_{1,n}(N)$ to an explicit but slightly messy problem in linear algebra.

The strategy can be outlined as follows. We have a forgetful mapping $\EuScript{M}_{1,n}(N)\to\EuScript{M}_{1,1}(N)$: we take a curve and forget all the marked points but one. (It does not matter which point we retain since the action of an elliptic curve on itself by translations is transitive.) The corresponding Leray spectral sequence degenerates in the second term for dimension reasons (the base is non-compact and 1-dimensional). Our first task will be describing this spectral sequence. Luckily, this is more or less standard, and the main techniques we'll be using are

\begin{itemize}
\item B. Totaro's interpretation \cite{totaro} of the Cohen-Taylor spectral sequence \cite{cohentayl} for the cohomology of the ordered configuration space of a smooth complex projective variety as the Leray spectral sequence of the open embedding of the configuration space into the cartesian power.
\item The ``correct'' way to introduce the (mixed) Hodge structure in the cohomology of a local system that
underlies a variation of (mixed) Hodge structures; we shall be interested only in the case when the base is a smooth curve, and we shall use the version
developed by
J. Steenbrink--S. Zucker \cite{sz} and S. Zucker \cite{zucker1}.
\item The Hodge-theoretic interpretation of the Eichler-Shimura isomorphism suggested by P. Deligne \cite{deligne} and proved by S. Zucker \cite{zucker}
(see also B\'ayer-NeuKirch \cite{baneu}).
\end{itemize}

In some cases minor modifications of these results will be needed. Thus, we'll show (theorem \ref{isoms}) that for an genus 1 curve $E$, the mixed Hodge structures on $H^*(F(E,n),\Q)$ coincide with those on the $E_\infty$ term of the Leray spectral sequence for $F(E,n)\subset E^{\times n}$ (i.e., the rational cohomology groups of the ordered configuration spaces of genus 1 curves decompose as direct sums of pure Hodge structures). Our proof works for any projective Abelian variety, and we conjecture that the result holds for any smooth complex projective variety $X$ with $b_1(X)\leq 2\dim X$.
We also compute (lemmas \ref{prrrop} and \ref{izo}) the rational mixed Hodge structure on the whole of the cohomology of the Hodge local systems on $\eu{M}_{1,1}(N)$, and not just on the image of the compactly supported cohomology (aka the Eichler cohomology).

As in \cite{getzler}, we first deal with the cohomology of the fine moduli spaces $\eu M_{1,n}(N),N>2$ by studying the Leray spectral sequence of the bundle $\eu M_{1,n}(N)\to\eu M_{1,1}(N)$ and then take the $\SL_2(\Z/N)$-invariant part.

The paper is organised as follows. In section \ref{secresults} we state our main results and sketch some proofs. In \ref{restop} we introduce some notation and describe the rational cohomology algebras of $\eu{M}_{1,n}$ (theorem \ref{prelim}) and $\eu{M}_{1,n}(N)$ in general (theorem \ref{main1}); the cup product determined up to a two term filtration in both cases.

The subsection \ref{Hodge} is devoted to the Hodge structures. In \ref{hhodge} we describe the complex mixed Hodge structures of $\eu{M}_{1,n}(N)$ (theorem \ref{main2}) and state a decomposition theorem for the rational cohomology of the configuration spaces of genus~1 curves (theorem \ref{isoms}). In \ref{hhodge} we also give an algorithm for computing the mixed Hodge polynomials of $\eu{M}_{1,n}(N)$ and interpret $H^*(F(E,n),\Q), E$ a genus~1 curve as the cohomology of a certain graph complex. In \ref{hodgelocsys} we describe the rational mixed Hodge structures on the cohomology of the Hodge local systems on $\eu{M}_{1,1}(N),N>2$ (theorem \ref{main3}) and on $\eu{M}_{1,1}(2)$ (corollary \ref{hlsm2}) and $\eu{M}_{1,1}(1)$ (corollary \ref{hlsm1}).

In the remaining sections we prove the results stated in \ref{secresults}.

I am grateful to Jozef Steenbrink for patiently answering many questions about the mixed Hodge theory and to Samuel Leli\`evre for a reference.

This is a preliminary version of the paper. The final version will be considerably shorter and include some examples of calculations (which this version does not). It will also include a discussion of the real homotopy type of the spaces $\eu{M}_{1,n}(N)$.

\section{Statement of results and discussion}\label{secresults}
\subsection{Rational cohomology of $\eu{M}_{1,n}$ as a topological space}\label{restop}

We begin with the topological part of the story. As promised above, our description will be quite complicated, so we'll need to introduce a lot of notation first.

\subsubsection{Cohomology of configuration spaces}\label{conf}
{\bf Notation.} For a topological space $X$, we denote by $F(X,n)$ the {\it space of ordered configurations of $n$ points on $X$}, i.e., the space formed by all $(x_1,\ldots,x_n)\in X^{\times n}$ such that $x_i\neq x_j$ for $i\neq j$.

Let $n$ be a positive integer. Let $\eu A_n$ be the graded commutative
$\Q$-algebra generated by the elements $\aaa_i,i\in\{1,\ldots,n\},\bb_i,i\in\{1,\ldots,n\}$ and $\triangle_{ij},i\neq j,i,j\in\{1,\ldots,n\}$ of degree 1 with the
following relations:

\begin{equation}\label{rel1}
\triangle_{ij}=\triangle_{ji},i,j\in\{1,\ldots,n\},i\neq j,\\
\end{equation}
\begin{equation}\label{rel2}
\triangle_{ij}\triangle_{jk}+\triangle_{jk}\triangle_{ki}+\triangle_{ki}\triangle_{ij}=0\mbox{ for any pairwise distinct indices }i,j,k\in\{1,\ldots,n\},
\end{equation}
$$\triangle_{ij}\aaa_i=\triangle_{ij}\aaa_j,\triangle_{ij}\bb_i=\triangle_{ij} \bb_j,i,j\in\{1,\ldots,n\},i\neq j.$$

We equip $\eu A_n$ with a differential graded algebra structure by introducing the differential $d$ defined by
$d(\aaa_i)=d(\bb_i)=0$ for any $i\in\{1,\ldots,n\}$ and $d(\triangle_{ij})=\bb_i\aaa_i+\bb_j\aaa_j-(\bb_i\aaa_j+\bb_j\aaa_i)$ for any $i,j\in\{1,\ldots,n\},i\neq j$. We shall
denote by $\eu H_n$ the cohomology algebra of $(\eu A_n,d)$. By the results of B. Totaro \cite{totaro}, $\eu H_n$ is the rational cohomology algebra of the space $F(E,n)$
of ordered configutations of $n$ points on an genus 1 curve $E$, and $\eu A_n$ is the $E_2$ term of the Leray spectral sequence of the open embedding $F(E,n)\subset E^{\times n}$. Note that $(\aaa_i,\bb_i)$ corresponds to a {\it negative} basis of the $H^1$ of the $i$-th factor in $E^{\times n}$; this choice will be explained later.

The elements $\sum_{i=1}^n\aaa_i,\sum_{i=1}^n\bb_i\in\eu A_n$ are obviously cocycles;
let $\eu I_n\subset \eu H_n$ be the ideal generated by the cohomology classes of these cocycles, and set $\eu B_n=\eu H_n/\eu I_n$.
An genus 1 curve $E$ acts on the
configuration space $F(E,n)$ by shifts. The classes in $\eu H_n$ of $\sum_{i=1}^n\aaa_i$ and $\sum_{i=1}^n\bb_i$ restricted to any orbit of this action span
the $H^1$ of the orbit. Moreover, the subalgebra of $\eu H_n$ generated by these classes is obvuously isomorphic to $H^*(E,\Q)$, and we get a splitting
\begin{equation}
H^*(F(E,n),\Q)\cong H^*(F(E,n)/E,\Q)\otimes H^*(E,\Q)
\end{equation}
that holds on the algebra level. The algebra
$\eu B_n$ can thus be naturally identified with $H^*(F(E,n)/E,\Q)$.

\subsubsection{The actions of $\SL_2(\Q)$ and $\mathfrak{S}_n$}\label{act}
Now we would like to define an action of $\SL_2(\Q)$ on $\eu A_n,\eu H_n$ and $\eu B_n$.
For an element $g=\bigl(\begin{smallmatrix}a&b&\\c&d\end{smallmatrix}\bigr)$ of $\SL_2(\Q)$ and $i,j\in\{1,\ldots,n\}$ we set
$$g\cdot \aaa_i=a\aaa_i+c\bb_i,g\cdot\bb_i=b\aaa_i+d\bb_i,g\cdot\triangle_{ij}=\triangle_{ij}.$$ It is easy to check that
this gives a well-defined action on $\eu A_n$ by algebra automorphisms; this action commutes with $d$. Hence, we obtain an action of $\SL_2(\Q)$ on $\eu H_n$ that obviously preserves $\eu I_n$ and hence,
descends to an action on $\eu B_n$.
%

If a group $G$ acts on a set $X$, we denote by $X^g$ the set of fixed points of $g\in G$; as usual, we set $X^G=\bigcap_{g\in G}X^g$. We shall write $X/G$ to denote the quotient space, regardless whether the action is left or right; we hope that this will never lead to a confusion.

The symmetric group $\mathfrak{S}_n$ acts on $\eu A_n$ by permuting the $\mathbf{a}_i$'s, the $\mathbf{b}_i$'s and the $\triangle_{ij}$'s. More precisely, for a $\sigma\in\mathfrak{S}_n$ we set $\sigma\cdot\mathbf{a}_i=\mathbf{a}_{\sigma(i)},\sigma\cdot\mathbf{b}_i=\mathbf{b}_{\sigma(i)},\sigma\cdot\triangle_{i,j}=\triangle_{\sigma(i),\sigma(j)}$. This action is easily seen to induce an action of $\mathfrak{S}_n$ on $\eu{H}_n$ and $\eu{B}_n$. Since the actions of $\mathfrak{S}_n$ and $\SL_2(\Q)$ on $\eu{B}_n$ commute, $\eu B_n^{\SL_2(\Q)}$ is a graded $\mathfrak{S}_n$-module. Now we can state our first result.

\begin{theorem}\label{prelim} 
The rational cohomology of $\eu{M}_{1,n}$ is additively isomorphic to $\eu{B}_n^{\SL_2(\Q)}\oplus\eu{D}^{*-1}_n$ with $\eu{D}^*_n$ the graded $\eu{B}_n$-module defined by
$$\eu{D}^*_n=\eu{B}_n^{\bigl(\begin{smallmatrix}-1&0&\\0&-1\end{smallmatrix}\bigr)}/\left(\eu{B}_n^{\bigl(\begin{smallmatrix}0&-1&\\1&0\end{smallmatrix}\bigr)}+\eu{B}_n^{\bigl(\begin{smallmatrix}0&-1&\\1&-1\end{smallmatrix}\bigr)}\right).$$ Under this identification, cup product of $(b,d)$ and $(b',d')$, $b,b'\in\eu{B}_n^{\SL_2(\Q)},d,d'\in \eu{D}^{*-1}_n$, is $$(bb',b\cdot d'+b'\cdot d+f_n(b,b')),$$ where $f_n$ is some 2-cocycle on $\eu{B}_n^{\SL_2(\Q)}$ (with coefficients in $\eu{D}^{*-1}_n$). The symmetric group $\mathfrak{S}_n$ acts on $\eu{D}^{*-1}_n$, and the $\mathfrak{S}_n$-action on $H^*(\eu{M}_{1,n},\Q)$ is obtained from those on $\eu{B}_n^{\SL_2(\Q)}$ and~$\eu{D}^{*-1}_n$.
\end{theorem}
%

In order to obtain later information about the mixed Hodge structures, we have to be a bit more specific about the actions of $\SL_2(\Q)$ on $\eu{A}_n,\eu{H}_n$ and $\eu{H}_n$ that we have just defined.
These actions are easily seen to be algebraic (meaning that the representation mappings from $\SL_2(\Q)$ to the groups of
$\Q$-automorphisms of the corresponding $\Q$-vector spaces are in fact algebraic). This implies that on each space $\eu A_n,\eu H_n,\eu B_n$
there is a natural action of the algebra $\mathfrak{sl}_2(\Q)$, and these spaces decompose as direct sums of the standard $\SL_2(\Q)$-modules
$\mathbf{V}_i$ \cite[Chapter VIII,\S 1]{bourbaki} (recall that $\mathbf{V}_i,i\geq 0$\label{V} is the $i$-th symmetric power of the standard 2-dimensional module). Note also that $\eu{A}_n^{\SL_2(\Q)}$ coincides with the $\Q$-subspace $\eu{A}_n^{\mathfrak{sl}_2(\Q)}$ consisting of the elements annihilated by all $A\in\mathfrak{sl}_2(\Q)$; the same holds for $\eu{H}_n$ and~$\eu{B}_n$.

Let $R_n:\SL_2(\Q)\to\Aut_\Q(\eu B_n)$ and $\rho_n:\mathfrak{sl}_2(\Q)\to\End_\Q(\eu B_n)$ be the representation morphisms.
Set $$\label{xyh}\mathbf{X}=\left(\begin{array}{cc}0&1\\0&0\end{array}\right), \mathbf{Y}=\left(\begin{array}{cc}0&0\\1&0\end{array}\right),
\mathbf{H}=\left(\begin{array}{cc}1&0\\0&-1\end{array}\right).$$

Given nonnegative integers $k$ and $i$, we define
$$\eu C^i_n(k)=\ker\rho_n(\mathbf{X})\cap\ker(\rho_n(\mathbf{H})-k\mathrm{Id}_{\eu B_n^i}),$$ where $\eu B_n^i$ is the degree $i$ component
of $\eu B_n$. This is, of course, the space of highest weight $k$ vectors of the representation $\rho_n$; i.e., $\eu C^i_n(k)$ ``counts''\label{counts} how many times the module $\mathbf{V}_k$ occurs in $\eu B_n^i$.
Notice that for any $k$ $$\eu C^*_n(k)=\bigoplus_i\eu C^i_n(k)$$ is both a $\eu B_n^{\SL_2(\Q)}$-module and a $\mathfrak{S}_n$-module.

\subsubsection{Spaces of automorphic forms}
As usual, we denote by $\Gamma(N),N>1$ the kernel of the natural morphism $\SL_2(\Z)\to\SL_2(\Z/N),$ and we set $\Gamma(1)=\SL_2(\Z)$. For a Fuchsian group $\Gamma$ of the first kind, we set $S_k(\Gamma)$, respectively $G_k(\Gamma)$, to be the space of weight $k$ cusp forms, respectively of weight $k$ modular forms, for~$\Gamma$. See e.g. \cite[chapter 2]{shimura} for the definition and explicit formulae for the dimensions of these spaces; note that modular forms are called ``integral forms'' in \cite{shimura}.
We set
$$g_{k,N}=\dim_\C G_k(\Gamma(N)),$$
$$s_{k,N}=\dim_\C S_k(\Gamma(N)).$$

Denote by $\mathbf{H}$ the upper half-plane $\{z\in\C\mid \Imm z>0\}$. In the sequel, when considering a subgroup of $\SL_2(\Z)$ or $\PSL_2(\Z)$ acting on $\mathbf{H}$, we shall mean the standard left action.


\subsubsection{The cohomology of $\eu{M}_{1,n}(N)$: the first version}

\begin{theorem}\label{main1} 
Let $f_n:\eu{M}_{1,n}(N)\to\eu{M}_{1,1}(N)$ be the forgetful morphism, and let $(E_r^{pq}(f_n))$ be the corresponding Leray spectral sequence. Then 

\begin{itemize} We have $E_2^{pq}(f_n)=0$ for $p\neq 0,1$, $E_2^{0,*}(f_n)\cong\eu{B}_n^{\SL_2(\Q)}$ as graded algebras, and
$E_r^{1,*}(f_n)$ is isomorphic to the graded $\eu{B}_n^{\SL_2(\Q)}$-module
\begin{equation}\label{e1*}
\bigoplus_{k\geq 0}(\eu C^*_n(k)\otimes\mathrm{W}(k,N))
\end{equation}
where $\mathrm{W}(k,N)$ is a $\Q$-vector space of dimension $s_{k+2,N}+g_{k+2,N}$ concentrated in degree 0 and equipped with the trivial action of $\eu{B}_n^{\SL_2(\Q)}$.

\item The cup-product on $(E_2^{pq}(f_n))$ in obtained from the product in $\eu{B}_n^{\SL_2(\Q)}$ and the action of $\eu{B}_n^{\SL_2(\Q)}$ on (\ref{e1*}).

\item The natural action of the symmetric group $\mathfrak{S}_n$ on $\eu{M}_{1,n}(N)$ covers the identical action on the base $\eu{M}_{1,1}(N)$ of $f_n$.
The resulting action of $\mathfrak{S}_n$ on $(E_2^{pq}(f_n))$ is obtained by using those on $\eu B_n^{\SL_2(\Q)}$ and $\eu C^*_n(k)$ described above and the trivial action on $\mathrm{W}(k,N)$.

\item Just for the record, let us state that the spectral sequence $(E_r^{pq}(f_n))$ degenerates in the second term i.e. $E_2(f_n)=E_3(f_n)=\cdots=E_\infty(f_n)$.
\end{itemize}
\end{theorem}

{\bf Remark.} The components $\mathrm{W}(k,N)$ are of course nothing else but the cohomology of $\Gamma(N)$ with coefficients in the standard $k+1$ dimensional module $\mathbf{V}_k$; in the case $N>2$ this is the same as the cohomology of the quotient of the upper half-plane by $\Gamma(N)$ with coefficients in the local system corresponding to $\mathbf{V}_k$.

The theorem implies that we have an exact sequence
\begin{equation}\label{extopo}
0\longrightarrow I\longrightarrow H^*(\eu{M}_{1,n}(N),\Q)\longrightarrow A\longrightarrow 0
\end{equation}
where $A\cong E_2^{0,*}(f_n)$ is a graded algebra and $I\cong E_2^{1,*-1}(f_n)$ is a graded $A$-module and an ideal in $H^*(\eu{M}_{1,1}(N),\Q)$. The general theory of spectral sequences tells us that the product
of any two elements from $I$ is zero, and the product $a'x$ of $x\in I$ and an element $a'$ lifting $a\in A$ is just $a\cdot x$. But a priori there is no reason the sequence (\ref{extopo}) should split multiplicatively, or equivalently, that $H^*(\eu{M}_{1,n}(N),\Q)$ contains a subalgebra that is mapped isomorphically onto $A$. We believe this is the case but for the moment let us state the following as a conjecture.

\newtheorem{conjecture}{Conjecture}

\begin{conjecture}\label{toposplit}
There exists a subalgebra of $H^*(\eu{M}_{1,n}(N),\Q)$ that restricts isomorphically onto $A\cong E^{0,*}\cong\eu B_n^{\SL_2(\Q)}$. Hence, $H^*(\eu{M}_{1,n}(N),\Q)\cong E_2(f_n)$ as algebras.
\end{conjecture}

\subsection{Mixed Hodge structures}\label{Hodge}
\subsubsection{Configuration spaces of elliptic curves}\label{hhodge}

Due to the above-mentioned geometric interpretation of $\eu A_n$, this algebra
is equipped with a mixed Hodge structure, which we now describe.

Let $\Lambda(n)$ be the graded commutative algebra generated by the degree 1 elements ${\triangle}_{ij},i,j=1,\ldots,n,i\neq j$ modulo relations (\ref{rel1}), (\ref{rel2}). This is in fact the rational cohomology algebra of $F(\C,n)$. We introduce the Tate Hodge structure of weight $2i$ on the degree $i$ part of $\Lambda(n)$. To a product $\triangle$ of the $\triangle_{ij}$'s there corresponds in an obvious way a partition $J(\triangle)$ of $\{1,\ldots,n\}$.
For such a partition, let $\Lambda_J(n)$ be the vector subspace of $\Lambda(n)$ spanned by the monomials corresponding to $J$, and let $\Lambda^q(n)$ be the degree $q$ part of $\Lambda(n)$; set $\Lambda^q_J(n)=\Lambda^q(n)\cap\Lambda_J(n)$. By definition, we have $\Lambda^q_J(n)=0$ if the number $|J|$ of the elements of $J$ is $<n-q$. On the other hand, it is easy to show (using $\Lambda^q(k)\cong H^q(F(\C,k),\Q)=0$ if $q>k-1$) that $\Lambda_J^q(n)=0$ if $|J|>n-q$.

Let us recall the isomorphism from \cite{totaro} between $\eu A_n$ and the $E_2$ term of the Leray spectral sequence for $F(E,n)\subset E^{\times n}$. 

{\bf Notation.}\label{subdiv} Let $J$ be a partition of $\{1,\ldots,n\}$ into $|J|$ subsets. The {\it multidiagonal} $E^{\times |J|}_J$ is the subvariety of $E^{\times n}$ defined by the condition that the ``coordinates'' with indices in the same element of the partition should be equal.

We have 
$$E^{pq}_2=\bigoplus_{|J|=n-q}\Lambda^q_J(n)\otimes H^p(E^{\times (n-q)}_J,\Q),$$ where $J$ runs through the set of partitions of $\{1,\ldots,n\}$ into $n-q$ subsets. Notice that $E_2^{pq}$ carries a natural pure Hodge structure of weight $p+2q$.

Let us describe the differential $d_2$. If $X$ is a smooth closed manifold, we have a natural mapping $H^*(X)\to H^{*+\dim_\R X}(X^{\times 2})$\label{xtox2} given as the composition of the mapping induced by either of the two projections and the cup-product with the class of the diagonal. Thus, if the partition $J'$ is obtained from $J$ by subdividing one of the elements into two non-empty subsets, we have a mapping $f_J^{J'}:H^*(E^{\times |J|}_J,\Q)\to H^{*+2}(E^{\times |J'|}_{J'},\Q)$ obtained from the obvious commutative square

$$
\begin{CD}
E^{\times |J|}_{J}@>>>E^{\times |J'|}_{J'}\\
@VVV @VVV\\
E^{\times |J|}@>{\mathrm{diag}\times\Id}>>E^{\times |J'|}
\end{CD}
$$

Now if $x=\triangle_1\cdots\triangle_q\otimes y\in E_2^{pq}$ where any $\triangle_l$ is one of the $\triangle_{ij}$'s, we have 
\begin{equation}\label{d2}
d_2(x)=\sum(-1)^{l+1}(\triangle_1\cdots\widehat{\triangle_l}\cdots\triangle_q)\otimes f_{J(\triangle_1\cdots\triangle_q)}^{J(\triangle_1\cdots\widehat{\triangle_l}\cdots\triangle_q)}(y).
\end{equation}

The cup product on $E_2$ is described as follows. Let $J_1$ and $J_2$, $|J_1|=n-i_1, |J_2|=n-i_2$ be two partitions
(where $J_1\cup J_2$ is the finest partition refined both by $J_1$ and $J_2$). We have the natural product mappings $$E^{\times (n-i_1-i_2)}_{J_1\cup J_2}\to E^{\times (n-i_1)}_{J_1}\times E^{\times (n-i_2)}_{J_2},$$ inducing the cohomology mappings
\begin{equation}\label{cupe2}
H^{p_1}(E^{\times (n-i_1)}_{J_1},\Q)\otimes H^{p_2}(E^{\times (n-i_2)}_{J_2},\Q)\to H^{p_1+p_2}(E^{\times (n-i_1-i_2)}_{J_1\cup J_2},\Q).
\end{equation}
On the other hand, the product in $\Lambda(n)$ gives us
\begin{equation}\label{prodlambd}
\Lambda_{J_1}(n)\otimes\Lambda_{J_2}(n)\subset\Lambda_{J_1\cup J_2}(n).
\end{equation}
It can be checked that $\Lambda_{J_1}(n)\cdot\Lambda_{J_2}(n)=0$ unless $|J_1\cup J_2|=n-i_1-i_2$.\footnote{Indeed, if $|J_1\cup J_2|\neq n-i_1-i_2$ , then $|J_1\cup J_2|> n-i_1-i_2$, so we have $\Lambda_{J_1}(n)\cdot\Lambda_{J_2}(n)=\Lambda_{J_1}^{i_1}(n)\cdot\Lambda_{J_2}^{i_2}(n)\subset\Lambda^{i_1+i_2}_{J_1\cup J_2}(n)=0$.} The tensor product of (\ref{prodlambd}) and (\ref{cupe2}) is equal to $(-1)^{p_1+q_2}$ times
the cup product $E_2^{p_1,q_1}\otimes E_2^{p_2,q_2}\to E_2^{p_1+p_2,q_1+q_2}$.

We have $E_2^{p,0}\cong H^p(E^{\times n},\Q)$, and $E_2^{0,1}$ is the sum of $\frac{n(n-1)}{2}$ 1-dimensional subspaces $\Lambda_J^1(n), J$ a partition that consists of $n-1$ elements.
An isomorphism $\eu A_n\to E_2$, which gives us the mixed Hodge structure on $\eu A_n$, is now constructed by picking a {\it negative} basis of $H^1(E,\Z)\subset H^1(E,\Q)$ and mapping the $\mathbf{a}_i$'s and the $\mathbf{b}_i$'s to $E_2^{1,0}$ and the $\triangle_{ij}$'s to $E_2^{0,1}$.

The differential $d$ on $\eu{A}_n$ is strictly compatible with the mixed Hodge structures that we have just described; we equip $\eu{H}_n$ and $\eu{B}_n$ with mixed Hodge structures induced from $\eu{A}_n$.

\medskip

\label{actx}
Let $e_1$ and $e_2$ be the elements of $H_1(E,\Z)$ that correspond to the elements $1$ and $\tau$ of the lattice $\langle 1,\tau\rangle$, and for $g=\bigl(\begin{smallmatrix}a&b&\\c&d\end{smallmatrix}\bigr)\in\SL_2(\Z)$ let $R_g$ be the $\R$-linear mapping $\C\to\C$ that takes 1 to $c\tau + d$ and $\tau$ to $a\tau + b$.
In this way we get a right action of $\SL_2(\Z)$ on $H_1(E,\Z)$
and a left action of $\SL_2(\Z)$ on $H^*(E,\Z)$. Let $(\aaa,\bb)$ be the basis of $H^*(E,\Z)$ such that $\aaa(e_1)=\bb(e_2)=0,\aaa(e_2)=\bb(e_1)=1$; an element $g\in\SL_2(\Z)$ as above takes $\aaa$ to $a\aaa+c\bb$ and $\bb$ to $c\aaa+d\bb$.

\begin{theorem}\label{isoms}
There exists an $\SL_2(\Z)\times\mathfrak{S}_n$-equivariant algebra isomorphism $\eu H_n\to H^*(F(E,n),\Q)$
that is compatible with the mixed Hodge structures and takes the classes of $\aaa_i$ and $\bb_i,i=1,\ldots,n$ to the restrictions to $F(E,n)$ of the pullbacks of $\aaa$ and $\bb$ under the $i$-th projection $E^{\times n}\to E$.
\end{theorem}

{\bf Remark.} Note that Totaro proves that $\eu H_n\cong H^*(F(E,n),\Q)$ as algebras using Deligne's grading of the weight filtration \cite[\S 4]{deligne1}. This grading in general not compatible with the Hodge filtration.
Hence, Totaro's argument does not allow one to compute the mixed Hodge numbers, and the mixed Hodge structure on $\eu{A}_n\cong E_2$ described above has apriori nothing to do with the one on $H^*(F(E,n),\Q)$. But luckily, the Leray spectral sequence of any morphism of complex algebraic varieties comes equipped from the second term on with a functorial mixed Hodge structure, which is compatible with all the differentials and the mixed Hodge structure on the cohomology of the source, see Saito \cite{saito} or Arapura \cite{arapura}. (Arapura considers only projective morphisms, but the result extends to the general case \cite[6.1]{stepet}.) Now, the functoriality alone easily shows that the mixed Hodge structure on the $E_2$ and $E_\infty$ the Leray spectral sequence of $F(E,n)\subset E^{\times n}$ 
coincides with the one described above.

It is easy to see that the Leray filtration on $H^*(F(E,n),\Q)$ that corresponds to the morphism $F(E,n)\subset E^{\times n}$ coincides up to a shift with the weight wiltration (cf. \cite[4]{totaro}), which implies that the action of $\SL_2(\Z)$ preserves the weight filtration on $H^*(F(E,n),\Q)$ and we have
$\eu H_n\cong \Gr^W_*(H^*(F(E,n),\Q))$ as algebras, $G$-modules and Hodge structures. So to prove theorem \ref{isoms} it would suffice to prove that $H^*(F(E,n),\Q)$ decomposes as a direct sum of pure Hodge substructures that are preserved under the action of $\SL_2(\Z)\times\mathfrak{S}_n$ (and the algebra isomorphism would follow automatically from the compatibility between the cup-product and the mixed Hodge structures). This will be done in subsection \ref{proofisoms}.


\medskip

Let $B_n$ be the image of the natural (algebra) monomorphism $H^*(F(E,n)/E,\Q)\to\eu{H}_n$; clearly, $B_n$ is a mixed Hodge substructure of $\eu{H}_n$ isomorphic to $H^*(F(E,n)/E,\Q)$.
The subalgebra $A$ of $\eu{H}_n$ generated by the classes of $\sum\aaa_i$ and $\sum\bb_i$ is a Hodge substructure too. We have $\eu{H}_n=\eu{I}_n\oplus B_n$ as mixed Hodge structures (since $\eu{I}_n$ is the image of $A^{\geq 1}\otimes B_n$ under the mixed Hodge structure isomorphism $A\otimes B_n\to\eu{H}_n$). We get the following corollary of theorem \ref{isoms}.
\begin{corollary} Equip $\eu{B}_n=\eu{H}_n/\eu{I}_n$ with a mixed Hodge structure induced from $\eu{H}_n$. We have then $\eu{B}_n\cong B_n\cong H^*(F(E,n)/E,\Q)$ as algebras, mixed Hodge structures and $\SL_2(\Z)$-modules.
\end{corollary}

$\clubsuit$

Notice that the weight filtrations on $\eu A_n, \eu H_n$ and $\eu B_n$ (as opposed to the Hodge filtrations) do not depend on the choice of $E$. Let us now take $E=\C/\langle 1,\mathrm{i}\rangle$. The Hodge filtrations will then be defined over $\Q(\mathrm{i})$.

\subsubsection{The forgetful morphism $\eu{M}_{1,n}(N)\to\eu{M}_{1,1}(N)$}
As noted above, the Leray spectral sequence of any morphism of complex algebraic varieties is equipped with a natural mixed Hodge structure.
We describe it here for the forgetful morphism $\eu{M}_{1,n}(N)\to\eu{M}_{1,1}(N)$

\begin{theorem}\label{main2} We use the notation of theorem \ref{main1}.
\begin{itemize}
\item The mixed Hodge structure on $E^{0,*}\cong\eu B_n^{\SL_2(\Q)}$ is induced by the inclusion $\eu B_n^{\SL_2(\Q)}\subset\eu B_n$. This mixed Hodge structure is the direct sum of Tate structures and does not depend on the choice of the genus 1 curve~$E$.
\item The mixed Hodge structure on $E^{1,*}\cong \bigoplus_{k\geq 0}(\eu C^*_n(k)\otimes\mathrm{W}(k,N))$ is induced by mixed Hodge structures on $\eu C^*_n(k)$ and $\mathrm{W}(k,N)$. The mixed Hodge structure on $\eu C^*_n(k)$ is described as follows.

Let $(\underline{W}_i)$ and $(\underline{F}^p)$ be the weight and the Hodge filtrations on $\eu B_n\otimes_\Q\Q(\mathrm{i})$. Identify $\eu C^{*}_n(k)$ with its image in $\eu B_n\otimes_\Q\Q(\mathrm{i})$ under the inclusion $c\mapsto c\otimes_\Q 1$. The mixed Hodge structure on $\eu C^{*}_n(k)$ is obtained by seting $$W_i=\underline{W}_{i-k}\cap\eu C^{*}_n(k),$$ $$F^p=(\underline{F}^p\cap\eu C^{*}_n(k))\otimes_\Q\C\subset\eu C^{*}_n(k)\otimes_\Q\C.$$ This mixed Hodge structure is also the direct sum of Tate ones.

\item We describe here the complex mixed Hodge structure on $\mathrm{W}(k,N)$; the rational one will be described later (theorem \ref{main3}).

The quotient $W_i(\mathrm{W}(k,N))/W_{i-1}(\mathrm{W}(k,N))$ is nonzero, iff $i=k+1,2k+2$. We have $\dim W_{k+1}(\mathrm{W}(k,N))/W_{k}(\mathrm{W}(k,N))=2s_{k+2,N}$. The Hodge structure on $$W_{2k+2}(\mathrm{W}(k,N))/W_{2k+1}(\mathrm{W}(k,N))$$ is Tate, and the Hodge decomposition of $W_{k+1}(\mathrm{W}(k,N))\otimes\C$ contains only $(0,k+1)$ and $(k+1,0)$ components.
\end{itemize}
\end{theorem}

As in the case of theorem \ref{main1}, the general theory (\cite{saito} or \cite{arapura}) tells us that the sequence (\ref{extopo}) is a sequence of mixed Hodge structures, but a priori there is no guarantee it splits. And again, it turns out that it actually does.

\begin{theorem}\label{hodgesplit}
There exists a Hodge substructure of $H^*(\eu{M}_{1,n}(N),\Q)$ that restricts isomorphically to $A\cong E^{0,*}\cong\eu B_n^{\SL_2(\Q)}$. Hence, $H^*(\eu{M}_{1,n}(N),\Q)\cong E_2(f_n)$ as algebras equipped with mixed Hodge structures.
\end{theorem}

%

For a complex algebraic variety $V$, we define the {\it mixed Hodge polynomial} of $V$ to be
$$P_{\mathrm{mHdg}}(V)=\sum_{n,p,q}t^nu^pv^q\dim_\C(\Gr^p_F\Gr_{p+q}^WH^n(V,\C)).$$
By setting in this expression $v=u$, respectively, $u=v=1$, we get the Poincar\'e-Serre polynomial, respectively, the Poincar\'e polynomial, of $V$.
The {\it mixed Hodge polynomial of $V$ with compact supports} (which we denote by $P_{\mathrm{mHdg},c}$) is obtained by replacing $H^n$ by $H^n_c$ in the definition of $P_{\mathrm{mHdg}}$;
by specialising $P_{\mathrm{mHdg},c}$ at $t=-1$ we get the Serre characteristic of~$V$. If there is a finite group $G$ acting on $V$ by automorphisms,
we define the {\it equivariant mixed Hodge polynomial} $P^G_{\mathrm{mHdg}}(V)$ in an obvious way. I.e., we replace the dimension $\dim_\C(\Gr^p_F\Gr_{p+q}^WH^n(V,\C))$ (which is the image of $\Gr^p_F\Gr_{p+q}^WH^n(V,\C)$ in the Grothendieck group of finite-dimensional vector spaces) by the image of $\Gr^p_F\Gr_{p+q}^WH^n(V,\C)$ in the Grothendieck group of finite-dimensional $G$-modules.

Theorem \ref{main2} gives us the following algorithm of computing the mixed Hodge polynomials of $\eu{M}_{1,n}(N)$. We denote by $R(G)$ the ring of (complex finite-dimensional algebraic) representations of an algebraic group $G$. Recall that above we have explicitly described the second term $(E_2^{pq})$ of the Leray spectral sequence of the embedding $F(E,n)\to E^{\times n}$ and introduced an action of $\SL_2(\Z)\times\mathfrak{S}_n$ on it; this action extends to an action of $\SL_2(\C)\times\mathfrak{S}_n$ on the complexification of $(E_2^{pq})$. For any $p,q$ set $g_{p,q}$ to be the image of $E^{pq}_3\otimes\C=E^{pq}_\infty\otimes\C$ in $R(\SL_2(\C)\times\mathfrak{S}_n)$; we have $g_{p,q}=\sum_{i=0}^\infty g_{p,q}^i[\mathbf{V}_i]$ where $g_{p,q}^i\in R(\mathfrak{S}_n)$. The expression $$\sum_{j=0}^\infty t^j \sum_{p+q=j}\sum_i(uv)^\frac{p+q-i}{2}g^i_{p,q}[\mathbf{V}_i]$$ is in fact a polynomial $\in R(\SL_2(\C)\times\mathfrak{S}_n)[t,u,v]$ and is divisible by $1+t[\mathbf{V}_1]+t^2$. Denote the quotient by $H$ and decompose $H$ as above: $H=\sum_i h_i[\mathbf{V}_i]$ with $h_i\in R(\mathfrak{S}_n)[t,u,v]$. The equivariant mixed Hodge polynomial $P^{\mathfrak{S}_n}_{\mathrm{mHdg}}(\eu{M}_{1,n}(N))$ is then equal to
$$h_0+t\sum_{i=0}^\infty[s_{i+2,N}(u^{i+1}+v^{i+1})+(g_{i+2,N}-s_{i+2,N})(uv)^{i+1}]h_i.$$

\medskip

Finally, let us note that there is an alternative way of looking at the algebra $\eu{A}_n=E_2$. Let $\eu{F}$ be the set of all forests (disjoint unions of trees) with the set of vertices equal to $\{1,\ldots,n\}$. (Thus we allow ``free'' vertices not adjacent to any edge.) Set $A=H^*(E,\Q)$. To each $F\in\eu{F}$ associate the (unordered) tensor product $V(F)=A^{\otimes\mbox{{\scriptsize the components of }}F}$. If the partitions of $\{1,\ldots,n\}$ corresponding to $F_1,F_2\in\eu{F}$ coincide, we have a natural identification $V(F_1)\cong V(F_2)$.

\begin{figure}\centering
\epsfbox{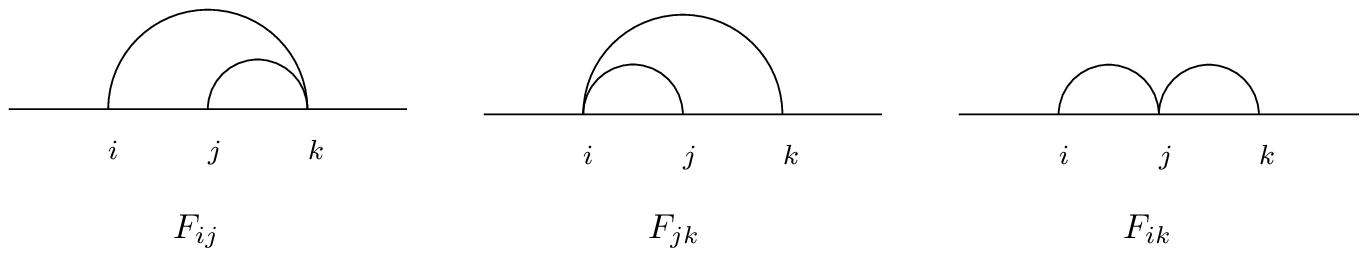}
\caption{}
\label{fig}
\end{figure}

The algebra $\eu{A}_n$ is additively isomorphic to $\bigoplus_{F\in\eu{F}}V(F)$ considered modulo the Arnold relation (\ref{rel2}), which can be translated in terms of forests as follows. Take an $F\in\eu{F}$, three indices $i,j,k$ whose the $F$-components are pairwise distinct, and let $F_{ij},F_{jk}$ and $F_{ki}$ be the graphs represented on Figure \ref{fig}; we set then $v_1+v_2+v_3=0$ for any $v_1\in V(F\cup F_{ij}),v_2\in V(F\cup F_{jk})$ and $v_3\in V(F\cup F_{ki})$ that are identified under $V(F_{ij})\cong V(F_{jk})\cong V(F_{ki})$. More precisely, to fix an isomorphism $\bigoplus_{F\in\eu{F}}V(F)/\mathrm{Arnold}\to\eu{A}_n\cong E_2$ let us assume that the subsets $\{i,j\}\subset\{1,\ldots,n\}$ are totally ordered, say, lexicographically. There is a natural isomorphism $\alpha_F:V(F)\cong H^*(E_J^{\times |J|},\Q)$ where $J$ is the partition corresponding to $F$, and we take $v\in V(F)$ to $\triangle\otimes\alpha_F(v)$ where $\triangle$ is the lexicographically ordered product of $\triangle_{ij}$'s that corresponds to $F$.
\captiondelim{. }

To see what becomes of the cup product, let us note that if $F=F_1\cup F_2, F, F_1, F_2\in\eu{F}$, then we have a natural map $\beta_{F_1,F_2}^F:V(F_1)\otimes V(F_2)\to V(F)$. Now, given arbitrary $F_1,F_2\in\eu{F}$, the product of the classes of $v_1\in V(F_1)$ and $v_2\in V(F_2)$ is zero if $F_1\cup F_2$ contains cycles and is otherwise equal to the class of $(-1)^{\varepsilon +|v_1|}\beta_{F_1,F_2}^{F_1\cup F_2}(v_1\otimes v_2)$ where $|v_1|$ is the degree of $v_1$ and $\varepsilon$ accounts for the difference between the lexicographical ordering of the edges of $F_1\cup F_2$ and another ordering obtained by first taking the ordered edges of $F_1$ and then the ordered edges of $F_2$.

Let us also translate the formula (\ref{d2}) for the differential. If $F'\in\eu{F}$ is obtained from $F\in\eu{F}$ by deleting one edge, we have a morphism $\gamma_F^{F'}:V(F)\to V(F')$ induced by the mapping $H^*(E,\Q)\to H^{*+2}(E^{\times 2},\Q)$ described on page \pageref{xtox2}. For the class of $v\in V(F), F\in\eu{F}$ we have then $$d(v)=\sum_{i}(-1)^{i+1}\gamma_F^{F_{l_i}}(v)$$ where $l_i$ is the $i$-th edge of $F$ and $F_{l_i}$ is $F$ minus $l_i$.

\subsubsection{Rational Hodge structures on $\mathrm{W}(k,N)$}\label{hodgelocsys}
Theorem \ref{main2} takes care of the complex mixed Hodge structures. The most explicit way to describe the rational ones is probably via group cohomology. Given a group $G$ and a $G$-module $V$, we define the cochains $C^*(G,V)$, the cocycles $Z^*(G,V)$ and the coboundaries $B^*(G,V)$ using the standard bar resolution (see, e.g., \cite{brown} or \cite{feiginfuchs}).

Consider a finite index torsion free subgroup $\Gamma\subset\PSL_2(\Z)$. A {\it special polygon} for $\Gamma$ (see \cite{kulkarni}\footnote{Special special polygons are defined in \cite{kulkarni} for an arbitrary finite index subgroup $\Gamma\subset\PSL_2(\Z)$; in the sequel we shall not need the general case.}) is a convex hyperbolic polygon $P\subset\mathbf{H}$ equipped with a fixed point free involution $\sigma$ of the set of the sides of $P$ such that
\begin{enumerate}
\item All vertices of $P$ belong to $\Q\cup\{\infty\}$.
\item $\Gamma$ is freely generated by the elements of $\PSL_2(\Z)$ obtained by picking, for any pair of sides of $P$ that are identified under $\sigma$, a transformation $\in\PSL_2(\Z)$ taking one to the other in an ``orientation-reversing'' way\footnote{The choice is unique, up to taking the inverse, since we are considering transformations from $\PSL_2(\Z)$.}, see \cite[section 2.5]{kulkarni} for details.
\end{enumerate}

The free generators of $\Gamma$ obtained from the second condition, and the inverses thereof, will be called {\it side-pairing transformations}.

For every $\Gamma$ as above there exists a special polygon \cite[theorem 3.3]{kulkarni}. Special polygons for the groups $\Gamma(N),N\leq 12$ were constructed in \cite{kulkarni} (we identify $\Gamma(N)$ with its image in $\PSL_2(\Z)$); in fact, \cite{kulkarni} describes some kind of a procedure that ``involves trial and error'' and enables one, given an integer $N$ (and a certain amount of luck), to construct a special polygon for $\Gamma(N)$.

Vertices of $P$ are, of course, cusps of $\Gamma$. If $v$ is a vertex of $P$, then a generator of the parabolic subgroup $\Gamma_v=\{g\in\Gamma\mid g\cdot v=v\}$ can be written as a product of side-pairing transformations as follows. Let $l$ be a side of $P$ adjacent to $v$, and let $g_1$ be the side pairing transformation that takes $l$ to $\sigma(l)$. If $g_1(v)\neq v$, then let $g_2$ be the side-pairing transformation that takes $l'$ to $\sigma(l')$, with $l'$ the side of $P$ adjacent to $g(v)$ and different from $g_1(l)$, etc. We continue in this way, until we have found the first $i$ such that $(g_i\cdots g_1)\cdot v=v$; the element $g_i\cdots g_1$ generates $\Gamma_v$.  

Let us fix an $N>2$. Take a special polygon $P,\sigma$ for $\Gamma(N)$, and construct a system $S=\{g_1,\ldots,g_s\}$ of free generators of $\Gamma(N)$ as above. Let $T:\Gamma(N)\to\Aut_\Q(\mathbf{V}_k)$ be the representation morphism for the module $\mathbf{V}_k$ (which was defined on page \pageref{V}). Let $(\mathbf{e}_1,\mathbf{e}_2)$ be the standard basis of $\mathbf{V}_1$.

\begin{theorem}\label{main3}
There exists an isomorphism \begin{equation}\label{iso}\mathrm{W}(k,N)\to H^1(\Gamma(N),\mathbf{V}_k)\end{equation} that satisfies the following.
\begin{enumerate}
\item The subspace $W_{k+1}(\mathrm{W}(k,N))$ is identified with $$\bigcap_v\ker(H^1(\Gamma(N),\mathbf{V}_k)\to H^1(\Gamma(N)_v,\mathbf{V}_k))$$ with $v$ running through the set of vertices of $P$.\footnote{One could as well let $v$ run through the set of $\Gamma(N)$-nonequivalent vertices of $P$.}
\item For a modular automorphic form $f$ of weight $k+2$ for $\Gamma(N)$ and any $z_0\in\mathbf{H}$, the function $\mathbf{S}_f:\Gamma(N)\to\mathbf{V}_k\otimes\C$ defined by $$\mathbf{S}_f(g)=\int_{z_0}^{g\cdot z_0}f(z)(z\mathbf{e}_1+\mathbf{e}_2)^kdz$$ is a cocycle $\in Z^1(\Gamma(N),\mathbf{V}_k\otimes\C)$. (The integration path can be chosen to be any rectifiable path in $\mathbf{H}$ joining $z_0$ and $g\cdot z_0$.)

If we complexify (\ref{iso}), the subspace $F_{k+1}(\mathrm{W}(k,N)\otimes\C)$, respectively, $$F_{k+1}(\mathrm{W}(k,N)\otimes\C)\cap W_{k+1}(\mathrm{W}(k,N))\otimes\C,$$ is identified with the subspace spanned by the classes of $\mathbf{S}_f, f\in G_{k+2}(\Gamma(N))$, respectively, by the classes of $\mathbf{S}_f, f\in S_{k+2}(\Gamma(N))$. 
\end{enumerate}
\end{theorem}

Explicitly, we can view $Z^1(\Gamma(N),\mathbf{V}_k)$ as the $\Q$-vector space of all functions $S\to\mathbf{V}_k$; indeed, since $\Gamma(N),N>2$ is free, each such function $\alpha$ can be extended to a cocycle $\alpha':\Gamma(N)\to\mathbf{V}_k$ by \begin{equation}\label{cocycle}
\alpha'(g_1g_2)=T(g_1)\alpha'(g_2)+\alpha'(g_1).\end{equation} The 1-coboundaries will correspond to the functions given by $g_i\mapsto T(g_i)x-x$ for some $x\in\mathbf{V}_k$. A function $\alpha:S\to\mathbf{V}_k$ gives an element of $\ker(H^1(\Gamma(N),\mathbf{V}_k)\to H^1(\Gamma(N)_v,\mathbf{V}_k))$, iff the corresponding cocycle $\alpha'$ satisfies $\alpha'(g_v)=T(g_v)x-x$ for some $x\in\mathbf{V}_k$ (here $g_v$ is a generator of $\Gamma(N)_v$.

Now let us consider the remaining cases $N=2,1$. We shall need the following construction. Let $G$ be a group, $H\subset G$ a subgroup, $V$ a $G$-module and $\res:H^*(G,V)\to H^*(H,V)$ the restriction morphism. If $(G:H)<\infty$, there exists a morphism $\tr^G_H:H^*(H,V)\to H^*(G,V)$ called the {\it transfer} that satisfies $\tr^G_H\circ\res=(G:H)\Id_{H^*(G,V)}$, see, e.g. \cite[chapter~III, \S 9]{brown}. This morphism is induced by the cochain morphism $\tr^G_H:C^*(H,V)\to C^*(G,V)$ constructed as follows (ibid.; we consider only the case of 1-cochains, since this is all we shall need).

Choose a representative $\bar g$ for any class $Hg, g\in G$, and define a mapping $\rho:G\to H$ by $\rho(g)=g\bar g^{-1}$. Let $\alpha:G\to V$ be a 1-cochain $\in C^1(H,V)$; we set $$\tr^G_H(\alpha)(g)=\sum_{\bar g}\bar g\alpha(\rho(\bar g^{-1} g)),$$ where the sum is taken over the representatives $\bar g$ that we have chosen.

Let $\bar\Gamma(2)$ be the image of $\Gamma(2)$ in $\PSL_2(\Z)$.
\begin{corollary}\label{hlsm2}
If $k$ is odd, then $\mathrm{W}(k,2)\cong H^1(\Gamma(2),\mathbf{V}_k)=0$.

Suppose $k$ is even. Equip $H^1(\Gamma(2l),\mathbf{V}_k),l>0$ with a mixed Hodge structure using theorem \ref{main3}. There exists an isomorphism
\begin{equation}\label{isog2}
\mathrm{W}(k,2)\to H^1(\bar\Gamma(2),\mathbf{V}_k)
\end{equation} such that the composition the mixed Hodge structure on $H^1(\bar\Gamma(2),\mathbf{V}_k)$ induced by (\ref{isog2}) coinsides with the one induced by $$\tr^{\bar\Gamma(2)}_{\Gamma(2l)}:H^1(\Gamma(2l),\mathbf{V}_k)\to H^1(\bar\Gamma(2),\mathbf{V}_k).$$ (We consider $\Gamma(2l)$ as a subgroup of $\PSL_2(\Z)$; note that $\bar\Gamma(2)$ acts on $\mathbf{V}_k$, since $k$ is even.)
\end{corollary}

Since $\bar\Gamma(2)$ is a free group, this gives an explicit description of the rational mixed Hodge structure on $\mathrm{W}(k,2)$.

Analogous statement holds in the case $N=1$.
\begin{corollary}\label{hlsm1}
Equip $H^1(\Gamma(l),\mathbf{V}_k),l>2$ with a mixed Hodge structure using theorem \ref{main3}. There exists an isomorphism
\begin{equation}\label{isosl2}
\mathrm{W}(k,1)\to H^1(\Gamma(1),\mathbf{V}_k)
\end{equation} such that the composition the mixed Hodge structure on $H^1(\Gamma(1),\mathbf{V}_k)$ induced by (\ref{isosl2}) coinsides with the one induced by $$\tr^{\Gamma(1)}_{\Gamma(l)}:H^1(\Gamma(l),\mathbf{V}_k)\to H^1(\Gamma(1),\mathbf{V}_k).$$
\end{corollary}

To make this completely explicit, let us describe the cohomology of $\Gamma(1)=\SL_2(\Z)$ with arbitrary rational coefficients.
Since $\SL_2(\Z)=\Z/4*_{\Z/2}\Z/6$, the Mayer-Vietoris sequence (see e.g., \cite[chapter VII, \S 9]{brown}) implies the following observation (which explains by the way theorem \ref{prelim}).

\begin{observation}
Let $V$ be a $\Q$-vector space equipped with an action $R:\SL_2(\Z)\to\Aut_\Q(V)$. Then
\begin{equation}\label{iso1}
H^1(\SL_2(\Z),V)\cong V^{ \bigl(\begin{smallmatrix}-1&0&\\0&-1\end{smallmatrix}\bigr)}/\left(V^{\bigl(\begin{smallmatrix}0&-1&\\1&0\end{smallmatrix}\bigr)}+V^{\bigl(\begin{smallmatrix}0&-1&\\1&-1\end{smallmatrix}\bigr)}\right).
\end{equation}
\end{observation}

This isomorphism can be easily constructed on the cocycle level. Namely, let $\alpha:\SL_2(\Z)\to V$ be a 1-cocycle, and set $$x_1=\alpha\bigl(\begin{smallmatrix}0&-1&\\1&0\end{smallmatrix}\bigr),x_2=\alpha\bigl(\begin{smallmatrix}0&-1&\\1&-1\end{smallmatrix}\bigr).$$
Let $A:V\oplus V\to V\oplus V$ be the mapping $$\left(R\bigl(\begin{smallmatrix}0&-1&\\1&0\end{smallmatrix}\bigr)-\Id_V\right)\oplus\left( R\bigl(\begin{smallmatrix}0&-1&\\1&-1\end{smallmatrix}\bigr)-\Id_V\right),$$ and choose a $\Q$-linear mapping $B:\Imm A\to V\oplus V$ such that $A\circ B=\Id_{\Imm A}$.
If $R\bigl(\begin{smallmatrix}0&-1&\\1&0\end{smallmatrix}\bigr)$ and $R\bigl(\begin{smallmatrix}0&-1&\\1&-1\end{smallmatrix}\bigr)$ are diagonalisable over $\C$ (which is the case for the modules $\mathbf{V}_k$), there is a distinguished choice for $B$. Since $\alpha|_{\Z/4}$ and $\alpha|_{\Z/6}$ are coboundaries, $x_1\oplus x_2\in\Imm A$. The element $B(x_1+x_2)$ is easily checked to be invariant under $R\bigl(\begin{smallmatrix}-1&0&\\0&-1\end{smallmatrix}\bigr)$, and the isomorphism (\ref{iso1}) is induced by the mapping that takes $\alpha$ to $B(x_1\oplus x_2)$.

\subsection{Unsettled questions and discussion}

\subsubsection{Mixed Hodge polynomials}

Theorems \ref{main1} and \ref{main2} give an algorithm to compute the equivariant mixed Hodge polynomials $P^{S_n}_{\mathrm{mHdg}}(\eu M_{1,n})$ and $P^{S_n}_{\mathrm{mHdg}}(\eu M_{1,n}(N))$. 
However, it would be interesting to obtain closed expressions or recurrent formulae for these polynomials. The main obstacle here is that while it is relatively easy to obtain the mixed Hodge polynomial of $\eu{A}_n=E_2$ of the Leray spactral sequence of $F(E,n)\subset E^{\times n}$, cf. \cite{getzler}, it is hard to see what happens when we pass to the cohomology.

\subsubsection{Relationship to the cell decomposition of the moduli spaces}

It would also be interesting to understand the relationship between the description of $H^*(\eu M_{1,n},\Q)$ given by theorems \ref{prelim} and \ref{main1} and the Harer-Mumford-Thurston-Penner cell decomposition (of the one-point compactification of $\eu M_{1,n}$), see \cite{harer,bowdepstein,penner}.

\subsubsection{Relationship to the tautological subalgebras}
It is well known (see e.g., \cite[exercises 3.28 and 4.9]{vakil}) that the tautological subalgebra of $H^*(\overline{\eu M}_{1,n},\Q)$ lives on the boundary $\overline{\eu M}_{1,n}\setminus\eu M_{1,n}$. Hence, the results of the present note are in a sence complementary (or maybe orthogonal\ldots) to those on the tautological subalgebras (see e.g., recent expositary papers \cite{morita} and \cite{vakil} for a description of the topological, respectively the geometric, approach).

\section{Proof of theorem \ref{isoms}}\label{proofisoms}
%


As remarked above, to prove theorem \ref{isoms}, it suffices to prove that $H^*(F(E,n),\Q)$ decomposes as a direct sum of pure Hodge substructures that are preserved under the action of $\SL_2(\Z)\times\mathfrak{S}_n$. By the Poincar\'e duality, it suffices to construct such a decomposition for $H^*(E^{\times n},D^n,\Q)$.

Let $D_1,\ldots D_{\binom{n}{2}}$ be the diagonals in $E^{\times n}$. We have a spectral sequence $(E_r^{pq})$ converging to $H^*(E^{\times n},D^n,\Q)$ with $E_1^{0,*}=H^*(E^{\times n},\Q)$,
$$E_1^{l,*}=\bigoplus_{1\leq j_1<\cdots< j_l\leq\frac{n(n-1)}{2}} H^*(D_{j_1}\cap\cdots\cap D_{j_l},\Q),l>0$$
and the differential $d_1$ given by $(d_1 \alpha)_j=\alpha|_{D_j}, \alpha\in E_1^{0,*}$,
$$(d_1\alpha)_{j_1,\ldots,j_{l+1}}=\left.\sum (-1)^{i+1}\alpha_{j_1,\ldots,\hat{j_i},\ldots,j_{l+1}}\right|_{D_{j_1}\cap\cdots\cap D_{j_{l+1}}}$$ for $$\alpha=\bigoplus_{1\leq j_1<\cdots< j_l\leq\frac{n(n-1)}{2}}\alpha_{j_1,\ldots,j_l}\in E_1^{i,*},i>0.$$ (See e.g. \cite{bengit}.)

All differentials in this spectral sequence are morphisms of mixed Hodge structures and the resulting filtration on $H^*(E^{\times n},D^n,\Q)$ is a filtration by mixed Hodge substructures. Since $E^{pq}_1$ is pure of weight $q$, the spectral sequence collapses at $E_2$.

Given an integer $l>1$, let $T_l$ be the endomorphism of $E$ defined by $a\mapsto l\cdot a$, and let $T_l^{\times n}$ be the corresponding diagonal automorphism of $E^{\times n}$. Since $T_l^{\times n}$ preserves all the diagonals, it induces an endomorphism (in fact, an automorphism) of $E_1$ compatible with $d_1$, and also an automorphism of $H^*(E^{\times n},D^n,\Q)$ which respects the filtration coming from the spectral sequence.
The endomorphism of $E_r^{pq}$ induced by $T_l^{\times n}$ is $l^q$ times identity.

Recall that if $\lambda$ is an eigenvalue of an endomorphism $f:V\to V$ of a finite-dimensional vector space $V$, then the generalised eigenspace corresponding to $\lambda$ is the set of all $v\in V$ which are annihilated by $(f-\lambda\mathrm{\mathop{Id}}_V)^{\dim V}$. The generalised eigenspaces of a morphism of rational mixed Hodge structure are mixed Hodge substructures.

Since $T_l^{\times n}$ commutes with $\SL_2(\Z)\times\mathfrak{S}_n$,
the decomposition of $H^*(E^{\times n},D^n,\Q)$ into the generalised eigenspaces of $(T_l^{\times n})^*$ splits the filtration coming from the spectral sequence both as a mixed Hodge structure and a $\SL_2(\Z)\times\mathfrak{S}_n$-module.$\clubsuit$


\section{Proofs of theorems \ref{main1} and \ref{main2}, $N>2$}\label{mainn>2}
In this section we assume $N>2$.

The moduli space $\eu{M}_{1,n}(N)$ is a locally trivial fibre bundle over $\eu{M}_{1,1}(N)\cong\mathbf{H}/\Gamma(N)$. For $\tau\in\mathbf{H},$ let $E_\tau$ denote the genus 1 curve $\C/\langle 1,\tau\rangle$. The $i$-th direct image of $\Q$ under this mapping is the local system with fibre $H^i(F(E_\tau,n)/E_\tau,\Q),\tau\in\mathbf{H}$. Let us compute the
monodromy of this bundle. Take a $g=\bigl(\begin{smallmatrix}a&b&\\c&d\end{smallmatrix}\bigr)\in\Gamma(N)$ and set $E'=E_\tau, E''=E_{g\cdot\tau}$.
We have a natural isomorphism $E'\to E''$ given by 

\begin{equation}\label{isoe1e2}
\C/\langle 1,\tau\rangle = \C/\langle c\tau+d,a\tau+b\rangle\to\C/\langle 1,\frac{a\tau+b}{c\tau+d}\rangle.
\end{equation}

Let $e_1$ and $e_2$ be the elements of $H_1(E',\Z)$ that correspond to the elements $1$ and $\tau$ of the lattice $\langle 1,\tau\rangle$. The monodromy mapping induced by the loop corresponding to $g$ takes $e_1$ to $d e_1+ce_2$ and $e_2$ to $be_1+ae_2$, i.e., we get a right action of $\Gamma(N)$ on $H_1(E',\Z)$ which is the composition of the inclusion $\Gamma(N)\subset\SL_2(\Z)$ and the action of $\SL_2(\Z)$ described above on page \pageref{actx}. (It is no surprise that the action is a right one, since the fundamental group of the base of a fibre bundle acts on the right on homotopy type of the fibre.) Recall that if $(\aaa,\bb)$ is the basis of $H^*(E',\Z)$ such that $\aaa(e_1)=\bb(e_2)=0,\aaa(e_2)=\bb(e_1)=1$, then the monodromy induced by $g$ takes $\aaa$ to $a\aaa+c\bb$ and $\bb$ to $b\aaa+d\bb$. Note that the subspace of holomorphic forms in $H^1(E',\C)$ is spanned by $\tau\aaa+\bb$.

%

Due to theorem \ref{isoms}, we know the mixed Hodge structure on
$H^i(F(E_\tau,n)/E_\tau,\Q)\cong\eu{B}_n$, as well as the algebra structure and the action of $\SL_2(\Z)$. A direct check shows that the natural action of $\SL_2(\Z)$ on $\eu{A}_n$ can be extended to an algebraic action of $\SL_2(\Q)$.

Let us now describe the mixed Hodge structure on $\eu{A}_n$ more explicitly.

Given an element $\tau\in\mathbf{H}$ and an integer $k\geq 0$, we define the rational mixed Hodge structure $\mathbf{V}_k^\tau$ of weight $k$ as follows. Let $(\mathbf{e}_1,\mathbf{e}_2)$ be the standard basis of $\Q^2$. The underlying $\Q$-vector space of $\mathbf{V}_k^\tau$ is set to be $\mathbf{V}_k=\Sym^k\Q^2$, and we define the $(p,q)$-component of $\mathbf{V}_k^\tau\otimes\C=\Sym^k\C^2$ as the one-dimensional subspace spanned by $$(\tau\mathbf{e}_1+\mathbf{e}_2)^p\cdot(\bar\tau\mathbf{e}_1+\mathbf{e}_2)^q.$$

As usual, for a $\Q$-vector space $V$ equipped with a mixed Hodge structure, we let $V(i),i\in\Z$ denote the Tate twist $V\otimes\Q(i)$ of $V$.

\begin{Prop}\label{mhsan}
The cohomology group $H^i(E_\tau^{\times k},\Q),i\geq 0,k\geq 1$ decomposes as a direct sum of $\mathbf{V}^\tau_{i-2j}(-j)$'s, $j=0,\ldots, \left[\frac{i}{2}\right]$.
\end{Prop}

This is an exercise in linear algebra, but since this statement will be important for the sequel, we shall give a proof.

{\bf Proof.}
The Hodge structure on $\mathbf{V}_k^\tau$ can be described in terms of representation theory: there exists a torus $T\subset\SL_2(\C)$ (that depends on $\tau$) and an isomorphism $f:\C^*\to T$ such that $f(z)\in\C^*$ acts on the $(p,q)$-component of $\mathbf{V}_k^\tau\otimes\C$ via multiplication by $z^{p-q}$.

Let us now take $\mathbf{V}^\tau_l$ and $\mathbf{V}^\tau_l,0\leq l\leq m$. To prove the proposition, it suffices to show that $\mathbf{V}^\tau_l\otimes\mathbf{V}^\tau_m,0\leq l\leq m$ decomposes as $$\bigoplus_{j=0}^{l}\mathbf{V}^\tau_{m-l+2j}(j-l)$$ as Hodge structures.
Write the weight decompositions of $\mathbf{V}^\tau_l$ and $\mathbf{V}^\tau_m$ with respect to $T$:
$$\mathbf{V}^\tau_l\otimes\C=V_{-l}\oplus\cdots\oplus V_l,$$
$$\mathbf{V}^\tau_m\otimes\C=W_{-m}\oplus\cdots\oplus W_m$$ where $f(z),z\in\C^*$ acts on $V_i$ and $W_i$ via multiplication by $z^i$.
We have \begin{equation}\label{hdg}
\mathbf{V}^\tau_l\otimes\mathbf{V}^\tau_m=\mathbf{V}_{l-m}\oplus \mathbf{V}_{l-m+2}\oplus\cdots\oplus \mathbf{V}_{l+m}
\end{equation} as $\SL_2(\Q)$-modules.

We have $$H^{p,q}(\mathbf{V}^\tau_l\otimes\mathbf{V}^\tau_m)=\bigoplus_{i+i=p-q}V_i\otimes W_j=\bigoplus_{j=0}^l\left(\mathbf{V}_{m-l+2j}^\C\cap\left( \bigoplus_{i+j=p-q}V_i\otimes W_j\right)\right)$$ with $\mathbf{V}_{m-l+2j}^\C=\mathbf{V}_{m-l+2j}\otimes\C$. But the space $\mathbf{V}_{m-l+2j}^\C\cap\left( \bigoplus_{i+j=p-q}V_i\otimes W_j\right)$ is the same as $$H^{\frac{{m-l+2j}+p-q}{2},\frac{{m-l+2j}+q-p}{2}}(\mathbf{V}^\tau_{m-l+2j})=H^{p,q}(\mathbf{V}^\tau_{m-l+2j}(j-l)),$$ and the proposition follows. $\clubsuit$

\begin{Prop}\label{mhshn}
Let $V$ and $W$ be direct sums of rational Tate Hodge structures equipped with the trivial actions of $\SL_2(\Q)$, and let $f:\mathbf{V}_k^\tau\otimes V\to \mathbf{V}^\tau_k\otimes W$ be a $\Q$-linear morphism of mixed Hodge structures. Then
the kernel and the image of $f$ are direct sums of Tate twists of $\mathbf{V}_k^\tau$. 
\end{Prop}

{\bf Proof.} This is another exercise in linear algebra. It suffices to show that if $\Q(i)$ and $\Q(j)$ are equipped with the trivial $\SL_2(\Q)$-actions, than any $\SL_2(\Q)$-equivariant morphism $g:\mathbf{V}_k^\tau\otimes\Q(i)\to\mathbf{V}_k^\tau\otimes\Q(j)$ of Hodge structures is zero unless $i=j$. If $i>j$, then this is true since the weight of $\mathbf{V}_k^\tau\otimes\Q(i)$ is smaller than the weight of $\mathbf{V}_k^\tau\otimes\Q(j)$. If $j>i$, then it follows from the compatibility of $g$ with the Hodge filtrations that $g$ has a nonzero kernel and hence is zero by the irreducibility of $\mathbf{V}_k$.
$\clubsuit$

Proposition \ref{mhsan} implies that as an $\SL_2(\Q)$-module equipped with a mixed Hodge structure, $\eu{A}_n$ decomposes as a direct sum of Tate twisted $\mathbf{V}^\tau_k$'s. By proposition \ref{mhshn}, so does $\eu{H}_n$ (the differential $d$ on $\eu{A}_n$ is strictly compatible with the mixed Hodge structures and it is easily checked that $d$ is $\SL_2(\Q)$-equivariant).

\begin{Prop}\label{mhsbn} Let $H$ be the direct sum $\bigoplus H^i$ of pure rational Hodge structures. Assume that the weights of the $H^i$'s are parewise distinct, and let $V$ be a rational mixed Hodge substructure of $H$. Then $V=\bigoplus (V\cap H^i)$ as Hodge structures.
\end{Prop}

{\bf Proof.} Writing the Deligne splittings, we see that $V\otimes\C=\bigoplus((V\otimes\C)\cap (H^i\otimes\C)).$ $\clubsuit$

Recall that $\eu{H}_n\cong\eu{I}_n\oplus\eu{B}_n$ as mixed Hodge structures. Since this splitting is compatible with the $\SL_2(\Q)$-action, we see (applying proposition \ref{mhsbn}) that $\eu{B}_n$ too decomposes as the direct sum of Tate twists of $\mathbf{V}_k^\tau$'s (both as a mixed Hodge structure and an $\SL_2(\Q)$-module)\label{decomppage}.

This implies that to calculate
the second term of the Leray spectral sequence of $\eu{M}_{1,n}(N)\to\eu{M}_{1,1}(N),N>2$, it remains to compute the cohomology of $\mathbf{H}/\Gamma(N)$ with coefficients in the polarised variation $\mathbb{V}_k$ of Hodge structure \cite{zucker} determined by the Hodge decompositions of $\mathbf{V}^\tau_k$. Notice that any $g\in\SL_2(\R)$ maps the $(p,q)$-compoment of $\mathbf{V}_k^\tau$ to the $(p,q)$-component of $\mathbf{V}_k^{g\cdot\tau}$, hence we can define the variation $\mathbb{V}_k$ for the quotient of the upper half-plane by any Fuchsian subgroup of the first kind without elliptic elements.

\begin{lemma}\label{prrrop} Let $\Gamma$ be a Fuchsian subgroup of $\SL_2(\R)$ of the first kind without elliptic points or irregular cusps. Assume that $-\Id\not\in\Gamma$. Then
\begin{enumerate}
\item The group $H^0(\mathbf{H}/\Gamma,\mathbb{V}_k)$ is $\Q(0)$ if $k=0$ and zero otherwise.
\item The subquotient $W_i(H^1(\mathbf{H}/\Gamma,\mathbb{V}_k))/W_{i-1}(H^1(\mathbf{H}/\Gamma,\mathbb{V}_k))$ is nonzero, iff $i=k+1,2k+2$.
\item The Hodge decomposition of the subspace $W_{k+1}(H^1(\mathbf{H}/\Gamma,\mathbb{V}_k))\otimes\C$ contains only $(k+1,0)$- and $(0,k+1)$-components, each of dimension $\dim S_{k+2}(\Gamma)$.
\item The subquotient $W_{2k+2}(H^1(\mathbf{H}/\Gamma,\mathbb{V}_k))/W_{2k+1}(H^1(\mathbf{H}/\Gamma,\mathbb{V}_k))$ is Tate of dimension equal to the number of cusps of $\Gamma$ for $k>0$ and to the number of cusps of $\Gamma$ minus 1 for $k=0$. (Recall that this is equal to $\dim_\C G_{k+2}(\Gamma)-\dim_\C S_{k+2}(\Gamma)$, see e.g. \cite[theorems 2.23-2.25]{shimura}.)
\item We have $H^{\geq 2}(\mathbf{H}/\Gamma,\mathbb{V}_k)=0$, if $\Gamma$ has a cusp or if $k>0$; if $\mathbf{H}/\Gamma$ is compact, then $H^2(\mathbf{H}/\Gamma,\mathbb{V}_0)\cong\Q(-1)$.
\end{enumerate}
\end{lemma}

{\bf Proof.} The first assertion of the lemma follows from the fact that there is no element of $\mathbf{H}$ fixed by all $g\in\Gamma$, and the last assertion is obvious.

Let $j$ be the embedding $\mathbf{H}/\Gamma\to S$ where $S$ is the natural compactification of $\mathbf{H}/\Gamma$ obtained by adding a point for each cusp of $\Gamma$ (for the details see e.g., \cite[\S 1.5]{shimura}). We will refer to the $\Gamma$-orbits of the cusps of $\Gamma$ as the cusps of $S$.
Assertions 2-4 of the lemma will be proved by studying the Leray spectral sequence for $j$.

The variation $\mathbb{V}_k$ is a direct summand of a geometric one, hence the properties (3.13) of \cite{sz} are satisfied by \cite[\S 5]{sz}, and Leray spectral sequence of $j$ is a spectral sequence of mixed Hodge structures \cite[theorem~4.1]{sz}.
It was proven by Zucker \cite[\S 12]{zucker} that $H^1(S,R^0 j_*\mathbb{V}_k)=H^1(S,j_*\mathbb{V}_k)$ is pure of weight $k+1$, the $(k+1,0)$ and $(0,k+1)$ components of $H^1(S,j_*\mathbb{V}_k)\otimes\C$ are of dimension $\dim S_{k+2}(\Gamma)$ and all other components vanish. This implies assertions 2-4 of the lemma in the case when the quotient $\mathbf{H}/\Gamma$ is compact. Assume from now on that $\Gamma$ has at least one cusp.

Let us now compute $H^0(S,R^1 j_*\mathbb{V}_k)$.

Set $\eu{V}_k=\eu{O}_{\mathbf{H}/\Gamma}\otimes\mathbf{V}_k$ and let $\bar{\eu V}_k$ be the canonical extension of $\eu{V}_k$ to $S$, see e.g., \cite[p. 91]{delequadiff}.
The complex $$K_*=[\bar{\eu V}_k\stackrel{\nabla}\to\bar{\eu V}_k\otimes\Omega^1_S(\log\Sigma)]$$ represents $Rj_*\mathbf{V}_k\otimes\C$ in the derived category $D^+(S)$ of sheaves on $S$ (here $\nabla$ is the canonical extension of the Gauss-Manin connection $\Id\otimes d$ and $\Sigma$ is the union of the cusps of $S$). This complex is equipped with the structure of a mixed Hodge complex of sheaves in the following way (see \cite[\S 4]{sz}): we set $F^p(K_*)=[\bar{\eu F}^p\stackrel{\nabla}\to\bar{\eu F}^{p-1}\otimes\Omega^1_S(\log\Sigma)]$, and
$$
W_i(K_*)=\begin{cases}
0, & i<k,\\
K_*' & k\leq i\leq 2k,\\
K_* & i>2k.
\end{cases}
$$
(Here $\bar{\eu F}^p$ the intersection of $\bar{\eu V}_k$ and the direct image of the Hodge bundle $\eu{F}^p$ on $\mathbf{H}/\Gamma$ and $K'_*$ denotes the complex $[\bar{\eu V}_k\to\Imm\nabla]$.)

This induces a mixed Hodge structure on $H^*(\mathbf{H}/\Gamma,\mathbb{V}_k)$, and also on the cohomology of $R^0j_*\mathbb{V}_k\otimes\C\stackrel{\mathrm{qis}}\sim\ker\nabla\stackrel{\mathrm{qis}}\sim K'_*$ and $R^1j_*\mathbb{V}_k\otimes\C\stackrel{\mathrm{qis}}\sim\coker\nabla\stackrel{\mathrm{qis}}\sim K_*/K'_*[1]$. The complex $K_*/K'_*$ is concentrated at the cusps of $S$. A simple topological argument shows that for a cusp $s$ of $S$ the group $H^0(s,R^1j_*\mathbb{V}_k)$ is one-dimensional; we claim the following

{\bf Claim.} The group $H^0(s,R^1j_*\mathbb{V}_k)$ is Tate of weight $2k+2$.

To prove this, it would suffice to show that for any $s\in\Sigma$ there exists a neighbourhood $U$ of $s$ and a section $\sigma$ of ${\eu F}^k$ over $U\setminus\{s\}$ of $s$ that extends to an section of $\bar{\eu V}_k$ over $U$ and such that $\frac{dt}{t}\otimes \sigma\not\in\Imm\nabla$ where $t$ is a local parameter of $S$ at $s$. Moreover, it suffices to look at what is happening at the cusp $s_0$ that corresponds to the infinity. Indeed\label{indeed}, if $s$ is any cusp of $\Gamma$, and $g\in\SL_2(\R)$ takes $s$ to the infinity, then the mapping $z\mapsto g\cdot z$ induces an isomorphism $\mathbf{H}/\Gamma\to\mathbf{H}/\Gamma'$ with $\Gamma'=g\Gamma g^{-1}$ that extends to an isomorphism between the natural compactifications of $\mathbf{H}/\Gamma$ and $\mathbf{H}/\Gamma'$. By considering the diagonal action of $g$ on $\Sym^k(\C^2)\times\mathbf{H}$ we obtain an isomorphism $\beta:{\eu V}_k\to\eu{V}_k'$ that covers $\alpha$ and extends to an isomorphism $\bar{\eu V}_k\to\bar{\eu V}_k'$ where $\bar{\eu V}_k'$ is the canonical extension of the pushforward ${\eu V}_k$ of $\Sym^k(\C^2)\otimes\eu{O}_\mathbf{H}$ to $\mathbf{H}/\Gamma'$. The isomorphism $\beta$ takes all Hodge subsheaves $\bar{\eu F}^p$ to the corresponding Hodge subsheaves.

Let us recall how to construct a local basis for the canonical extension of a flat holomorphic vector bundle.

Set $U=\{z\in\C\mid |z|<1\}$ and $U^*=U\setminus\{0\}$. If $\eu{V}$ is a holomorphic vector bundle over $U^*$ equipped with a flat connection with a unipotent monodromy, then there is a way to construct a local basis over ${\eu O}(U)$ of the canonical extension $\bar{\eu V}$ of $\eu V$ to $U$ starting from a basis of the space of flat multivalued sections of $\eu V$ (i.e., flat sections of the pullback of $\eu V$ to the universal cover of $U^*$), see e.g., \cite[proposition 11.3]{sz}. Namely, if $f:z\mapsto\mathrm{e}^{\frac{2\pi\mathrm{i}z}{h}},z\in\mathbf{H},h>0$ is the universal cover and $\mathbf{H}\ni u\mapsto \phi(u)$ is a flat multivalued section, then the corresponding section $\sigma_\phi$\label{secti} of $\bar{\eu V}$ is $f(u)\mapsto \mathrm{e}^{\frac{-u L}{h}}\phi(u)$ where $L$ is the (fibrewise) nilpotent logarithm of the counter-clockwise monodromy around zero. Let us apply this to our situation.

We can take $U^*=\mathbf{H}_r/\Gamma_0$ as a deleted neighbourhood of $s_0$, where $\mathbf{H}_r=\{ z\in\C\mid \Imm z>r\geq 0\}$ and $\Gamma_0\subset\Gamma$ is the subgroup generated by a matrix of the form $\bigl(\begin{smallmatrix}1&h&\\0&1\end{smallmatrix}\bigr),h>0$. (Here we use the fact that $\Gamma$ does not have irregular cusps.)
The universal covering map will then be just the projection $\mathbf{H}_r\to\mathbf{H}_r/\Gamma_0$. Take $\eu V$ to be the restriction of $\eu V_k$ to $U^*$. The pullback of $\eu V$ to $\mathbf{H}_r$ will be the bundle $\eu E=(\Sym^k(\C^2)\otimes{\eu O}(\mathbf{H}_r))$, and the fibrewise monodromy will act on each fibre $\eu E$ via $R_k\bigl(\begin{smallmatrix}1&-h&\\0&1\end{smallmatrix}\bigr)$ (recall that $R_k$ is the representation morphism associated to the module $\mathbf{V}_k\otimes\C$). If we apply the above procedure to the flat section $\tau\mapsto\mathbf{e}_2^k$ of $\eu E$, we get
the section $\tau\mapsto (\tau\mathbf{e}_1+\mathbf{e}_2)^k$ of $\eu E$ (cf. \cite[lemma 12.10]{zucker}), which pushes down to a section $\sigma$ of the Hodge bundle ${\eu F}^k$ restricted to $U^*$;\footnote{It is in fact a stroke of luck that we obtain precisely a section of the ${\eu F}^k$ in this way, see \cite[remark 12.13]{zucker}.} by the construction, $\sigma$ extends to a section of~$\bar{\eu V}_k$.

Now set	$U=U^*\cup\{s_0\}$ and take the section $\alpha=\frac{dt}{t}\otimes\sigma\in\Gamma(\bar{\eu V}_k\otimes\Omega^1_S(\log\Sigma),U)$ where $t=\mathrm{e}^{\frac{2\pi\mathrm{i}\tau}{h}}$ is a local parameter of $S$ at $s_0$. We have $\alpha\not\in\Imm\nabla$. Indeed, let $\mathbb{V}^{\lim}$ be the space of flat multivalued sections of ${\eu V}_k$ over $U^*$, and let $L$ be the nilpotent logarithm of the fibrewise monodromy around $s_0$ acting on ${\eu V}_k|U^*$. This extends to an endomorphism of $\bar{\eu V}_k|U$, see e.g., \cite[theorem II 1.17]{delequadiff}.
The Poincar\'e residue map takes $\Imm\nabla$ to $\Imm L$, see e.g., \cite[proof of proposition 11.3]{stepet}. On the other hand, the mapping $\mathbb{V}^{\lim}\ni\phi\mapsto\sigma_\phi\in\Gamma(\bar{\eu V}_k,U)$ is $L$-equivariant, so the image of $\alpha$ under the Poincar\'e map is
not contained in $\Imm L$ (since the section $\tau\mapsto\mathbf{e}_2^k$ is not). This proves the above claim about $H^0(s,R^1j_*\mathbb{V}_k)$.

To summarise, so far we have computed the terms $E_2^{0,1}$ and $E_2^{1,0}$ of the Leray spectral sequence for $j$; to complete the proof of assertions 2-4 of lemma \ref{prrrop}, we have to know the $E_2^{0,2}$ term, which is isomorphic to $H^2(S,j_*\mathbb{V}_k)\cong H^2_c(\mathbf{H}/\Gamma,\mathbb{V}_k)$. This group is zero for $k>0$.

If $k=0$, we have $H^2_c(\mathbf{H}/\Gamma,\mathbb{V}_k)=\Q$, and the kernel of the differential $d_2:E_2^{0,1}\to E_2^{2,0}$ is of codimension 1 (recall that we are considering the case when $\mathbf{H}/\Gamma$ is noncompact).

The lemma is proven.$\clubsuit$

{\bf Remark.}\label{Rjweight} It follows from the proof of lemma \ref{prrrop} that $W_{k+1}(H^1(\mathbf{H}/\Gamma,\mathbb{V}_k))\otimes\C$
coincides with the image of $H^1(S,R^0j_*\mathbb{V}_k\otimes\C)$ in $H^1(\mathbf{H}/\Gamma,\mathbb{V}_k\otimes\C).$

Propositions \ref{mhsan}, \ref{mhshn}, \ref{mhsbn} and lemma \ref{prrrop} allow one to write down the Leray spectral sequence for the bundle $\eu{M}_{1,n}(N)\to\eu{M}_{1,1}(N),N>2$ (this spectral sequence contains only the zeroth and the first columns). Thus, to prove
theorems \ref{main1} and \ref{main2} for $N>2$, it remains to explain the mixed Hodge structure on the spaces $\eu{C}^*_n(k)$.

Recall
that $\eu{C}^*_n(k)$ counts the multiplicity of $\mathbf{V}_k$ in $\eu B_n$ by counting the highest weight vectors with respect to the standard generators $\mathbf{X},\mathbf{Y},\mathbf{H}$ of $\mathfrak{sl}_2(\Q)$, see page \pageref{xyh}. The vector space $\eu{B}_n$ decomposes as a direct sum of Tate twisted $\mathbf{V}_k^\tau$'s, and the mixed Hodge structure on $\eu{C}^*_n(k)$ simply counts the twists; indeed, an easy check shows that the image in $\mathbf{V}_k^\tau\otimes K$ of an element of $\ker\mathbf{X}\subset \mathbf{V}_k^\tau(l)$ has Hodge filtration $-l$ (here $K$ is a subfield of $\C$ containing~$\tau$ and $\bar\tau$).

\section{Proof of theorem \ref{main3}}
In the previous section we proved that the spaces $F^{k+1}H^1(\mathbf{H}/\Gamma(N),\mathbb{V}_k\otimes\C)$ and $F^{k+1}W_{k+1}H^1(\mathbf{H}/\Gamma(N),\mathbb{V}_k\otimes\C)$ are isomorphic to $G_{k+2}(\Gamma(N))$ and $S_{k+2}(\Gamma(N))$ respectively. Here we show that these identifications can be made $\SL_2(\Z/N)$-equivariant.

\begin{lemma}\label{izo} Let $\Gamma$ be a subgroup of $\SL_2(\R)$ as in lemma \ref{prrrop}.Then there exists an isomorphism $G_{k+2}(\Gamma)\to F^{k+1}H^1(\mathbf{H}/\Gamma,\mathbb{V}_k\otimes\C)$ that identifies $S_{k+2}(\Gamma)$ with $F^{k+1}W_{k+1}H^1(\mathbf{H}/\Gamma,\mathbb{V}_k\otimes\C)$. If $\Gamma=\Gamma(N),N\geq 3$, this isomorphism can be chosen to be $\SL_2(\Z/N)$-equivariant.
\end{lemma}

{\bf Proof.} We use the notation from the proof of lemma \ref{prrrop}. The complex $F^{k+1}(K_*)$ is just $[0\to\bar{\eu F}^k\otimes\Omega^1_S(\log\Sigma)]$. We shall now describe the sections of $\bar{\eu F}^k\otimes\Omega^1_S(\log\Sigma)$.

To any modular form $f\in G_{k+2}(\Gamma)$ we associate the section
\begin{equation}\label{secc}
z\mapsto f(z)(z\mathbf{e}_1+\mathbf{e}_2)^kdz
\end{equation}
of the bundle $\Sym^k(\C^2)\otimes{\eu O}(\mathbf{H})\otimes\Omega^1(\mathbf{H})$. This section pushes down to a section $T(f)$ of ${\eu V}_k\otimes\Omega_{\mathbf{H}/\Gamma}^1$.
Indeed, under the natural left action of $\Gamma$ on $$\Gamma[\Sym^k(\C^2)\otimes{\eu O}(\mathbf{H})\otimes\Omega^1(\mathbf{H})]$$ the section (\ref{secc}) is transformed by an element $g^{-1}\in\Gamma$ into
\begin{equation}\label{transform}
z\mapsto f(g\cdot z)((g\cdot z)g^{-1}\cdot\mathbf{e}_1+g^{-1}\cdot\mathbf{e}_2)^kg^*(dz),
\end{equation}
which is equal to $f(z)(z\mathbf{e}_1+\mathbf{e}_2)^kdz.$
Moreover, at each $\tau\in\mathbf{H}$ the value of (\ref{secc}) belongs to the $(k,0)$-component of $\mathbf{V}_k^\tau$, i.e., $T(f)\in\Gamma({\eu F}^k\otimes\Omega^1_{\mathbf{H}/\Gamma})$.

Now let us show that $T(f)$ extends to a section of $\bar{\eu V}_k\otimes\Omega_{S}^1(\log\Sigma)$. Consider first the cusp $s_0$ at the infinity; as in the proof of lemma \ref{prrrop}, choose a generator $\bigl(\begin{smallmatrix}1&h&\\0&1\end{smallmatrix}\bigr),h>0$ of the corresponding parabolic subgroup of $\Gamma$ and a local coodrinate $u=\mathrm{e}^{\frac{2\pi\mathrm{i}}{h}z},\Imm z>r\geq 0$ for $S$ at $s_0$. As we have seen above, there exists a neighbourhood $U\ni s_0$ and a trivialisation $\bar{\eu V}_k|_U\cong\Sym^k(\C^2)\otimes\eu{O}_U$ that transforms the section $z\mapsto(z\mathbf{e}_1+\mathbf{e}_2)^k$ to the constant section $u\mapsto \mathbf{e}_2^k$. Tensoring by $\Omega_U^1(\log\{s_0\})$ we see that $T(f)$ locally is
\begin{equation}\label{locf}
f(\frac{h}{2\pi\mathrm{i}}\log u)\mathbf{e}_2^kdz=f(\frac{h}{2\pi\mathrm{i}}\log u)\mathbf{e}_2^k\frac{h}{2\pi\mathrm{i}}\frac{du}{u}.
\end{equation}
Since $f$ is a modular form, $u\mapsto f(\frac{h}{2\pi\mathrm{i}}\log u)$ is holomorphic at $0$, and $T(f)$ extends to a section of $\bar{\eu V}_k|_U\otimes\Omega_U^1(\log\{s_0\})$. Note that if $f$ is a cusp form, then this extension is in fact in $\bar{\eu V}_k|_U\otimes\Omega_U^1$.

Let now $s$ be any cusp of $\Gamma$, $g\in\SL_2(\R)$ an element taking $s$ to the infinity. Indeed as above (see p. \pageref{indeed}), set $\Gamma'=g\Gamma g^{-1}$, let $S'$ denote the natural compactification of $\mathbf{H}/\Gamma'$, and set $\bar{\eu V}_k'$ to be the canonical extension of the pushforward of $\Sym^k(\C^2)\otimes\eu{O}_\mathbf{H}$ to $\mathbf{H}/\Gamma'$. Let $\Sigma'$ be the set of the cusps of $S'$ and write $g^{-1}=\bigl(\begin{smallmatrix}a&b&\\c&d\end{smallmatrix}\bigr)$. If $f\in G_{k+2}(\Gamma)$, then $h:z\mapsto (cz+d)^{-k-2}f(g^{-1}\cdot z)$ is in $G_{k+2}(\Gamma')$ (see e.g., \cite[proposition 2.4]{shimura}), and the isomorphism $${\eu V}_k\otimes\Omega_{S}^1(\log\Sigma)\to {\eu V}'_k\otimes\Omega_{S'}^1(\log\Sigma')$$ that covers the natural isomorphism $\mathbf{H}/\Gamma\to\mathbf{H}/\Gamma'$ takes the section $T(f)$ to $T(h)$. Since $T(h)$ extends to a section of $\bar{\eu V}'_k\otimes\Omega_{S'}^1(\log\Sigma')$ at the infinity, $T(f)$ extends to a section of $\bar{\eu V}_k\otimes\Omega_{S}^1(\log\Sigma)$ at $s$.

To summarise, we have shown that $T(f)$ extends to a section of $\bar{\eu F}^k\otimes\Omega_S(\log\Sigma)$, which is in fact a section of $\bar{\eu F}^k\otimes\Omega_S$ if $f\in S_{k+2}(\Gamma)$. It is easy to check that the converse is also true: if $s\in\Gamma(\bar{\eu F}^k\otimes\Omega_S(\log\Sigma))$, then $s=T(f)$ for some $f\in G_{k+2}(\Gamma)$. Indeed, restricting $s$ to $\mathbf{H}/\Gamma$ and then pulling pack to $\mathbf{H}$ we obtain a section of $\eu{O}_\mathbf{H}\otimes\Omega^1_\mathbf{H}$ that can be written as $\bar s:z\mapsto f(z)(z\mathbf{e}_1+\mathbf{e}_2)^kdz$ for some holomorphic function $f\in{\eu O}_\mathbf{H}$. Since $\bar s$ is $\Gamma$-invariant, we see that $f$ satisfies $f(\gamma\cdot z)=(cz+d)^{k+2}f(z)$ for any $\gamma=\bigl(\begin{smallmatrix}a&b&\\c&d\end{smallmatrix}\bigr)\in\Gamma$, cf. (\ref{transform}). Moreover, since $s$ is a section of $\bar{\eu V}_k\otimes\Omega_S^1(\log\Sigma)$, we see that $f$ is holomorphic at the cusp at the infinity (cf. (\ref{locf})) and, in fact, at any cusp of $S$, i.e., $f\in G_{k+2}(\Gamma)$.

Thus $f\mapsto T(f)$ gives an isomorphism between $G_{k+2}(\Gamma)$ and $\Gamma(\bar{\eu F}^k\otimes\Omega_S(\log\Sigma))$. It follows from the construction of $T$ that if $\Gamma=\Gamma(N)$, then $T$ is $\SL_2(\Z/N)$-equivariant.
We know by lemma \ref{prrrop} that the dimension of $F^{k+1}H^1(S,\mathbb{V}_k\otimes\C)$ is equal to $\dim G_{k+2}(\Gamma)$, hence, the mapping
\begin{multline}\label{mapp}
G_{k+2}(\Gamma)\to\Gamma(\bar{\eu F}^k\otimes\Omega_S(\log\Sigma))=H^0(S,\bar{\eu F}^k\otimes\Omega_S(\log\Sigma))\to F^{k+1}\mathbb{H}^1(K_*)\\ \cong F^{k+1}H^1(\mathbf{H}/\Gamma,\mathbb{V}_k\otimes\C)
\end{multline}
that we have constructed
is an isomorphism (here $K_*$ is the complex of sheaves defined in the proof of lemma \ref{prrrop}). This isomorphism is $\SL_2(\Z/N)$-equivariant if $\Gamma=\Gamma(N)$. (To see this, recall that the isomorphism $\mathbb{H}^1(K_*)\cong H^1(\mathbf{H}/\Gamma,\mathbb{V}_k\otimes\C)$ is the composition
$$
\mathbb{H}^1(K_*)\to \mathbb{H}^1([j_*\eu{V}_k\to j_*\Omega_{\mathbf{H}/\Gamma}])\to\mathbb{H}^1([\eu{V}_k\to \Omega_{\mathbf{H}/\Gamma}])\gets H^1(\mathbf{H}/\Gamma,\mathbb{V}_k\otimes\C)
$$
where the first arrow is induced by the inclusion $K_*\subset [j_*\eu{V}_k\to j_*\Omega_{\mathbf{H}/\Gamma}]$, the second one is the restriction to $\mathbf{H}/\Gamma$ and the third one is induced by the morphism of $\mathbb{V}_k\otimes\C$ into its holomorphic de Rham resolution.)
 
Since we know (lemma \ref{prrrop}) that the $\dim F^{k+1}W_{k+1}H^1(\mathbf{H}/\Gamma,\mathbb{V}_k\otimes\C)=\dim_{k+2}(\Gamma)$, to complete the proof of lemma \ref{izo}, it remains to show that (\ref{mapp}) takes $S_{k+2}(\Gamma)$ to $F^{k+1}W_{k+1}H^1(\mathbf{H}/\Gamma,\mathbb{V}_k\otimes\C)$. This is essentially done in \cite[\S 12]{zucker} and \cite[\S 4]{baneu}. But since we work in a slightly different setting, let us give a proof.

Take a cusp form $f$; as we have seen above, we have then $T(f)\in\Gamma(\bar{\eu F}^k\otimes\Omega_S)$.
This implies that $T(f)\in\Imm\nabla$. Indeed, let $U$ be a neighbourhood of a cusp $s$ isomorphic to the unit disc, and let $\sigma_\phi$ be the section of $\bar{\eu V}_k|_U$ constructed from a flat multivalued section of $\eu V_k|_U$, see p. \pageref{secti}. Then (see e.g., \cite[proof of proposition 11.3]{stepet})
$$\nabla\sigma_\phi(u)=-N\phi(u)\otimes\frac{du}{u},u\in U^*=U\setminus\{s\}.$$ Applying the Leibnitz formula, we see that the question whether a section $\sigma\in\Gamma(\bar{\eu V}_k|_U\otimes\Omega^1_U(\log\{s\}))$ is in the image of $\nabla$ is reduced to the existence of a holomorphic solution $x:W\to\C^{k+1}$ of the equation
\begin{equation}\label{eqdiff}
\frac{dx}{du}=\frac{Ax}{u}+y(u)
\end{equation}
where $A$ is a nilpotent matrix of the maximal rank, $W\subset U$ is a neighbourhood of $s$ and $y$ is a meromorphic function on $U^*$. If $\sigma$ is a section of $\bar{\eu V}_k|_U\otimes\Omega^1_U$, then the corresponding function $y$ is holomorphic, and (\ref{eqdiff}) has a holomorphic solution.

So, the hypercohomology class $\in \mathbb{H}^1(K_*)$ corresponding to $T(f),f\in S_{k+2}(\Gamma)$ is in the image of $\mathbb{H}^1(F^{k+1}(K_*)\cap K'_*)$; this class has thus weight filtration $k+1$ (since $K'_*=W_kK_*$). Lemma \ref{izo} is proven.$\clubsuit$

Now we would like to identify $H^1(\mathbf{H}/\Gamma(N),\mathbb{V}_k\otimes\C)$ (computed using the complex of $\mathbf{V}_k\otimes\C$-valued differential forms on $\mathbf{H}/\Gamma(N)$) with $H^1(\Gamma(N),\mathbf{V}_k\otimes\C)$ (computed using the bar resolution, see e.g., \cite[I.5 and III.1]{brown}) in an $\SL_2(\Z/N)$-equivariant way. This can be done as follows.

\begin{Prop}\label{sl2nequi}
Take a base point $z_0\in\mathbf{H}$ and let $\omega=\sum_{i=0}^k \ee_1^i\ee_2^{k-i}\otimes\omega_i$ be a closed 1-form on $\mathbf{H}$ with values in $\mathbf{V}_k\otimes\C$. We put $\int_a^b\omega=\sum_{i=0}^k\left(\int_a^b\omega_i\right)\ee_1^i\ee_2^{k-i},a,b\in\mathbf{H}$. For any $g\in\SL_2(\R)$ set
\begin{equation}\label{actonforms}
g\cdot\omega=\sum_{i=0}^k (g\cdot\ee_1^i)(g\cdot\ee_2^{k-i})\otimes(g^{-1})^*\omega_i.
\end{equation}
If $\omega$ is invariant with respect to any element of $\Gamma(N)$, then the function $\mathbf{S}_\omega:\Gamma(N)\to\mathbf{V}_k\otimes\C$ defined by $\mathbf{S}_\omega(g)=\int_{z_0}^{g\cdot z_0}\omega$ is a {1-cocycle}; if moreover $\omega=df$ for a $\Gamma(N)$-invariant $C^{\infty}$-function $f:\mathbf{H}\to\mathbf{V}_k\otimes\C$, then $\mathbf{S}_\omega$ is a coboundary. The mapping $\omega\mapsto\mathbf{S}_\omega$ gives an $\SL_2(\Z/N)$-equivariant isomorphism $$\mathbf{S}^*:H^1(\mathbf{H}/\Gamma(N),\mathbb{V}_k\otimes\C)\to H^1(\Gamma(N),\mathbf{V}_k\otimes\C)$$ defined over $\Q$.
\end{Prop}

{\bf Proof.} An immediate check shows that $\mathbf{S}_\omega$ is a cocycle for $\omega$ closed, and that $\mathbf{S}_{df}$ is a coboundary, if $f\in C^{\infty}(\mathbf{H},\mathbf{V}_k\otimes\C)^{\Gamma(N)}$. To show that $\mathbf{S}^*$ is an isomorphism, let us note that any $\Gamma(N)$-invariant mapping from the orbit of $z_0$ to $\mathbf{V}_k\otimes\C$ can be extended to a smooth $\Gamma(N)$-invariant function $\mathbf{H}\to\mathbf{V}_k\otimes\C$ (and hence, any 1-coboundary of $\Gamma(N)$ with values in $\mathbf{V}_k\otimes\C$ is in fact $\mathbf{S}_{df}$ for some $f\in C^{\infty}(\mathbf{H},\mathbf{V}_k\otimes\C)^{\Gamma(N)}$, which implies that $\mathbf{S}^*$ is injective).

The mapping $\mathbf{S}^*$ is manifestly defined over $\Q$ and it remains to prove that $\mathbf{S}^*$ is $\SL_2(\Z/N)$-equivariant. Denote the space of $\C$-valued smooth 1-forms on $\mathbf{H}$ by $\eu{E}^1_\mathbf{H}$. The action of $g'\in\SL_2(\Z/N)$ on $\omega\in(\mathbf{V}_k\otimes{\eu E}^1_\mathbf{H})^{\Gamma(N)}$ is given by (\ref{actonforms}) with $g$ being an element of the preimage of $g'$ in $\SL_2(\Z)$. Let us now recall how $\SL_2(\Z/N)$ acts on $H^1(\Gamma(N),\mathbf{V}_k)$ (see e.g., \cite[III.8]{brown}).

Let $F_*$ be the standard resolution of $\Z$ over the group ring $\Z\Gamma(N)$ (see e.g., \cite[I.5]{brown}). For any $g\in\SL_2(\Z)$ the mapping $\alpha_g:F_*\to F_*$ given (in the standard $\Z$-basis) by $$\alpha_g(\gamma_0,\ldots,\gamma_l)=(g\gamma_0 g^{-1},\ldots g\gamma_l g^{-1})$$ is compatible (in the sence of \cite[II.6]{brown}) with the automorphism $\gamma\mapsto g\gamma g^{-1}$ of $\Gamma(N)$. The group $\SL_2(\Z)$ acts on the left on $\Hom(F_*,\mathbf{V}_k)^{\Gamma(N)}$ by $(g\cdot f)(x)=g\cdot f(\alpha_{g^{-1}}x),x\in F_*$. Identifying as usual $\Hom(F_1,\mathbf{V}_k)^{\Gamma(N)}$ with $\Map(\Gamma(N),\mathbf{V}_k)$ via $f\mapsto (\gamma\mapsto f(1,\gamma))$ we get an action of $\SL_2(\Z)$ on 
$C^1(\Gamma(N),\mathbf{V}_k)$ given by $(g\cdot F)(\gamma)=g\cdot F(g^{-1}\gamma g),F:\Gamma(N)\to\mathbf{V}_k$; this induces the required action of $\SL_2(\Z/N)$ on $H^1(\Gamma(N),\mathbf{V}_k)$.

Notice that for any $a,b\in\mathbf{H},g\in\SL_2(\R)$ and a closed 1-form $\omega\in\mathbf{V}_k\otimes\eu{E}^1_\mathbf{H}$ we have $$\int_a^bg\cdot\omega=g\cdot\int_{g^{-1}\cdot a}^{g^{-1}\cdot b}\omega.$$
Let us now take a closed form $\omega\in(\mathbf{V}_k\otimes\eu{E}^1_\mathbf{H})^{\Gamma(N)}$ and elements $g\in\SL_2(\Z)$ and $\gamma\in\Gamma(N)$. We have $$\mathbf{S}_{g\cdot\omega}(\gamma)=\int_{z_0}^{g\cdot z_0} g\cdot\omega=g\cdot \int_{g^{-1}\cdot z_0}^{g^{-1}\gamma\cdot z_0}\omega=g\cdot\left(\int_{z_0}^{g^{-1}\gamma g\cdot z_0}\omega+\int_{g^{-1}\gamma g\cdot z_0}^{g^{-1}\gamma\cdot z_0}\omega-\int_{z_0}^{g^{-1}\cdot z_0}\omega\right)$$ since $\omega$ is closed. The function $\gamma\mapsto g\cdot\int_{z_0}^{g^{-1}\gamma g\cdot z_0}\omega=g\cdot\left( \mathbf{S}_\omega(g^{-1}\gamma g)\right)$ is $g\cdot\mathbf{S}_\omega$ and the function
\begin{multline*}
\gamma\mapsto g\cdot\left(\int_{g^{-1}\gamma g\cdot z_0}^{g^{-1}\gamma\cdot z_0}\omega-\int_{z_0}^{g^{-1}\cdot z_0}\omega\right)=g\cdot\left(g^{-1}\gamma g\cdot\int_{z_0}^{g^{-1}\cdot z_0}\omega-\int_{z_0}^{g^{-1}\cdot z_0}\omega\right)\\=\gamma\left(g\cdot\int_{z_0}^{g^{-1}\cdot z_0}\omega\right)-g\cdot\int_{z_0}^{g^{-1}\cdot z_0}\omega
\end{multline*}
is a coboundary. This shows that the cohomology mapping $\mathbf{S}^*$ is $\SL_2(\Z/N)$-equivariant. Proposition \ref{sl2nequi} is proven.$\clubsuit$

It follows from the remark on page \pageref{Rjweight} and \cite[proposition 12.5]{zucker} that when we identify $H^1(\mathbf{H}/\Gamma(N),\mathbb{V}_k)$ with
$H^1(\Gamma(N),\mathbb{V}_k)$ using proposition \ref{sl2nequi}, the subspace
$W_{k+1}(H^1(\mathbf{H}/\Gamma(N),\mathbb{V}_k))$ is identified with $$\bigcap_v\ker(H^1(\Gamma(N),\mathbf{V}_k)\to H^1(\Gamma(N)_v,\mathbf{V}_k))$$ where $v$ runs through the set of the orbits of the cusps of $\Gamma(N)$, and $\Gamma(N)_v$ is the stabiliser of the cusp $v$. This observation, together with proposition \ref{prrrop} and lemma \ref{izo}, implies theorem~\ref{main3}.

\section{Proofs of theorems \ref{main1} and \ref{main2}, $N=1,2$}

To prove the remaining cases of theorems \ref{main1} and \ref{main2}, let us note that if $N_1$ divides $N_2$, then the group $\Gamma(N_1)/\Gamma(N_2)$ acts on the moduli space $\eu{M}_{1,n}(N_2)$, the quotient being $\eu{M}_{1,n}(N_1)$. Let us describe this action. Let $g$ be a lifting of an element $[g]$ to $\Gamma(N_1)\subset\SL_2(\Z)$ and take a $\tau\in\mathbf{H}$.
The action of $[g]$ on $\eu{M}_{1,n}(N_2)$ in obtained from the action on $\eu{M}_{1,1}(N_2)$ and the isomorphisms (\ref{isoe1e2}) for all $\tau$.

From this it follows that the action preserves (any) decomposition of $R^ip_*\Q$ (where $p:\eu{M}_{1,n}(N_2)\to\eu{M}_{1,1}(N_2)$ is the natural projection) as a sum of Tate twisted $\mathbb{V}_k$'s (see page \pageref{decomppage}). To show this if suffices to consider the case $N_1=1$; we have then $\Gamma(N_1)/\Gamma(N_2)=\SL_2(\Z/N_2)$. The images of $\bigl(\begin{smallmatrix}0&1&\\-1&0\end{smallmatrix}\bigr)$ and $\bigl(\begin{smallmatrix}0&-1\\1&-1\end{smallmatrix}\bigr)$ in $\SL_2(\Z/N_2)$ generate $\SL_2(\Z/N_2)$ and preserve any decomposition $$R^ip_*\Q=\bigoplus\mu(i,k,l)\mathbb{V}_k(l),\mu(i,k,l)\in\Z_{\geq 0}$$ over $[\mathrm{i}]$ and $\left[\mathrm{e}^{\frac{2\pi\mathrm{i}}{3}}\right]$ respectively (here $[\tau]$ stands for the image of $\tau\in\mathbf{H}$ in $\eu{M}_{1,1}(N_2)=\mathbf{H}/\Gamma(N_2)$).

The action of $[g]\in\Gamma(N_1)/\Gamma(N_2)$ on $\mathbb{V}_1$
can be described as follows when converted from the natural right action (since the fibers of $\mathbb{V}_1$ are cohomology rather than homology groups) to the left one. Take a representative $g=\bigl(\begin{smallmatrix}a&b&\\c&d\end{smallmatrix}\bigr)$ of $[g]$ as above, take two genus 1 curves $E'=E_\tau,\tau\in\mathbf{H}$ and $E''=E_{g\cdot\tau}$, let $e'_1,e'_2$ be the elements of $H_1(E',\Z)$ that correspond respectively to the elements 1 and $\tau$ of the lattice $\langle 1,\tau\rangle$, and define similarly $e''_1,e''_2\in H_1(E'',\Z)$. Finally, set $(\aaa',\bb')$ and  $(\aaa'',\bb'')$ to be the bases respectively of $H^1(E',\Z)$ and $H^1(E'',\Z)$ that are dual respectively to $(e'_2,e'_1)$ and $(e''_2,e''_1)$.

The mapping induced by $[g]$ takes $(e'_1,e_2')$ to $(ae_1''-ce_2'',-be_1''+de_2'')$, and hence, $(\aaa',\bb')$ is taken to $(a\aaa''+c\bb'',b\aaa''+d\bb'')$. The action of $\Gamma(N_1)/\Gamma(N_2)$ on all $\mathbb{V}_k$ is obtained from the action on $\mathbb{V}_1$ in the obvious way. A $\mathbb{V}_k\otimes\C$-valued form $T(f)$ on $\mathbf{H}/\Gamma(N_2)$ obtained from a modular form $f\in G_{k+2}(\Gamma(N_2))$ by (\ref{secc}) is transformed under $[g^{-1}]$ by formula (\ref{transform}). Theorems \ref{main1} and \ref{main2} for $N=1,2$ follow now from the fact that the space of $\Gamma(N_1)/\Gamma(N_2)$-invariant modular forms $\in G_{k+2}(\Gamma(N_2))$ coincides with $G_{k+2}(\Gamma(N_1))$, and analogously for cusp forms, see e.g., \cite[proposition 2.6]{shimura}.

\section{Proof of theorem \ref{hodgesplit}}

In this section we prove theorem \ref{hodgesplit}
We will use the relative version of the multiplication by $l$ endomorphism $T_l$, which was already used in section \ref{proofisoms} to prove that the cohomology of the configuration spaces of an elliptic curve decomposes as a direct sum of pure Hodge structures.

Let $\eu{E}(N)$  be the universal elliptic curve with a level $N$ structure. As an analytic variety, $\eu{E}(N)$ is the quotient of $\C\times\mathbf{H}$ by the action of the group generated by $(z,\tau)\mapsto (z+k+l\tau), k,l\in\Z$ and $(z,\tau)\mapsto \left(\frac{z}{c\tau+d},\frac{a\tau+b}{c\tau+d}\right)$ where $\left(\begin{smallmatrix}a&b&\\c&d\end{smallmatrix}\right)\in\Gamma(N)$. Note that this group is an extension of $\Gamma(N)$ by $\Z\oplus\Z$. There is a natural map $\eu{E}(N)\to X(N)=\mathbf{H}/\Gamma(N)$; let $\eu{F}(\eu{E}(N),n)$ be the fibered configuration space i.e. the subvariety of $F(\eu{E}(N),n)$ formed by the $n$-tuples $(x_1,\ldots,x_n)$ such that all $x_i$ lie in the same fibre. Set $\eu{E}^n(N)$ to be the $n$-th fibered power of $\eu{E}(N)$ over $X(N)$ and denote the fibered fat diagonal $\eu{E}^n(N)\setminus \eu{F}(\eu{E}(N),n)$ by $D^n(N)$. For the rest of the section fix a level $N>2$.

Let $E$ be a genus 1 curve. We will first prove theorem \ref{hodgesplit} by showing that the mixed Hodge substructure of $\SL_2(\Z)$-invariants $H^*(F(E,n),\Q)^{\SL_2(\Z)}$ extends to a mixed Hodge substructure of $H^*(\eu{F}(\eu{E}(N),n),\Q)$. 

Let $\eu{T}_l$ be the endomorphism of $\eu{E}^n(N)$ whose restriction to any fibre is multiplication by an integer $l\geq 2$ (cf. section \ref{proofisoms}). Rather than working with cohomology, we will work with the relative compactly supported cohomology groups $H^*_c(\eu{F}(\eu{E}(N),n))\cong H^*_c(\eu{E}^n(N),D^n(N))$ because $\eu{T}_l$ induces an endomorphism (in fact, automorphism) of these groups.

We denote the natural projection from $\eu{E}^n(N)$ to $X(N)$ by $p$.

\subsection{Proof of theorem \ref{hodgesplit}}

Consider the Leray spectral sequence with compact supports $(E^{i,j}_k,d_k)$ constructed from the map $p:\eu{E}^n(N)\to X(N)$ and the sheaf $$\eu{F}=j_!{j}^{-1}\underline{\Q}_{\eu{F}(\eu{E}(N),n)}$$ where $j$ is the embedding of $\eu{F}(\eu{E}(N),n)\subset\eu{E}^n(N)$. Note that $\eu{F}$ is isomorphic to the complex $$\left[ \underline{\Q}_{\eu{E}^n(N)}\to i_*\underline{\Q}_{D^n(N)}\right].$$

The map $\eu{T}_l$ acts on the spectral sequence in a natural way; the resulting endomorphism will be denoted $\eu{T}^*_l$. The higher derived images with compact supports of $\eu{F}$ under $\eu{T}_l$ are direct sums of Tate twisted $\mathbb{V}_k$'s extended to the cusp points by zero. The action of $\eu{T}^*_l$ on the Leray spectral sequence preserves this decomposition. Note that $p$ covers the identity of the base $X(N)$, so the action of $\eu{T}_l$ on the compactly supported cohomology of the sheaves $R^jp_!\eu{F}$ is determined by the action on the stalks.
If $\eu{G}\cong\mathbb{V}_k(-i)$ is a direct summand of $R^j p_! \eu{F}$, then $\eu{T}_l^*$ acts on $\eu{G}$ by multiplication by $l^{2i+k}$ (regardless of $j$). Note also that the Leray spectral sequence has only two nonzero columns, the first and the second, since $H^m_c(X(N),\mathbb{V})=0$ unless $m=1$ or 2.

By the Poincar\'e duality and the calculation of the cohomology of $\mathbb{V}_k$, see lemma \ref{prrrop}, the weight filtration on $H^1_c(X(N),\mathbb{V}_k)$ jumps at $0$ and $k+1$, i.e. if $W_i=W_iH^1_c(X(N),\mathbb{V}_k)$, then $W_i/W_{i-1}\neq 0$ iff $i=0,k+1$. So we can represent $H^1_c(X(N),\mathbb{V}_k(-i))$ as a direct sum of two pure Hodge structures $V_1\oplus V_2$ such that $T_l^*$ acts on $V_1$ by multiplication by $l^{w_1-1}$ and on $V_2$ by $l^{w_2+k}$ where $w_m,m=1,2$ is the weight of a nonzero element of $V_m$.

If $k\neq 0$, then $H^2_c(X(N),\mathbb{V}_k)=0$. The group $H^2_c(X(N),\mathbb{V}_0)\cong \Q(-1)$ is pure of weight~2. So $T_l^*$ acts on $H^2_c(X(N),\mathbb{V}_0(i))$ by multiplication by $l^{w-2}$ where $w$ is the weight of a nonzero element of $H^2_c(X(N),\mathbb{V}_0(i))$.

Let $E_\alpha\subset H^*_c(\eu{F}(\eu{E},n),\Q)$ be the generalised eigenspace of $T_l^*$ corresponding to the eigenvalue $l^\alpha$. Each $E_\alpha$ is a mixed Hodge substructure.

\medskip

{\bf Claim.} $\bigoplus_\alpha W_{\alpha+1} E_\alpha$ projects isomorphically onto the first column of the Leray spectral sequence of $p$.

{\bf Proof of the claim.}
Set $E^i_\alpha$ to be the degree $i$ part of $E_\alpha$. The Leray filtration gives the exact sequence of mixed Hodge structures

$$0\to E_\alpha^i\cap E^{2,i-2}_2\to E_\alpha^i \to \mbox{ the image of }E^i_\alpha\mbox{ in } E^{1,i-1}_2\to 0.$$

To simplify the notation let us temporarily denote $\underline{\Q}_{\eu{F}(\eu{E}(N),n)}$ as $\underline{\Q}$. The intersection $E_\alpha^i\cap E^{2,i-2}_2$ is the generalised eigenspace of $E^{2,i-2}_2$ corresponding to $l^\alpha$. This is a Hodge substructure isomorphic to
$$\bigoplus_{\begin{array}{c}\mathbb{V}_0(-j)\subset R^{i-2}p_!\underline{\Q} \\ 2j=\alpha\end{array}}H^2_c(X(N),\mathbb{V}_0(-j)),$$
which is pure of weight $2+2j=2+\alpha$.

The image of $E^i_\alpha$ in $E^{1,i-1}_2$ is the generalised eigenspace of $E^{1,i-1}_2$ corresponding to $l^\alpha$, which is isomorphic to
$$\bigoplus_{\begin{array}{c}\mathbb{V}_k(-j)\subset R^{i-1}p_!\underline{\Q} \\ 2j+k=\alpha\end{array}}H^1_c(X(N),\mathbb{V}_k(-j)).$$
A non-zero element of the latter group has weight $2j=\alpha-k$ or $2j+k+1=\alpha+1$; in both cases the weight is $<\alpha+2$.

Since morphisms of mixed Hodge structures are strict, the weight subspace $W_{\alpha+1}E^i_\alpha$ projects onto the image of $E^i_\alpha$ in $E^{1,i-1}_2$. Since a non-zero element of $E_\alpha^i\cap E^{2,i-2}_2$ has weight $\alpha+2$, the intersection $W_{\alpha+1}E^i_\alpha\cap(E_\alpha^i\cap E^{2,i-2}_2)=0$ so the map from $W_{\alpha+1}E^i_\alpha$ to the image of $E^i_\alpha$ in $E^{1,i-1}_2$ is injective. So $W_{\alpha+1}E_\alpha^i$ projects isomorphically to the generalised eigenspace of $E^{1,i-1}_2$ corresponding to $l^\alpha$. The claim is proven.

Theorem \ref{hodgesplit} follows from the claim since the Leray spectral sequence for cohomology is dual to the Leray spectral sequence for compactly supported cohomology.

\begin{flushright}
{\sc Alexei Gorinov\\
Department of Mathematics\\
National Research University ``Higher School of Economics''\\
Moscow\\
Russia}\\
\url{agorinov@hse.ru}\\
\url{gorinov@mccme.ru}
\end{flushright}
\end{document}